\documentclass{article}
\usepackage{amssymb}
\usepackage{amsfonts}
\usepackage{amsmath}

\newtheorem{corollary}{Corollary}[section]

\newtheorem{definition}{Definition}[section]

\newtheorem{lemma}{Lemma}[section]

\newtheorem{proposition}{Proposition}[section]

\newtheorem{theorem}{Theorem}[section]
\newenvironment{proof}[1][Proof]{\noindent\textbf{#1.} }{\ \rule{0.5em}{0.5em}}

\input{tcilatex}
\begin{document}

\title{On MF Property of Reduced Amalgamated Free Products of UHF Algebras }
\author{Qihui Li \ \ \ \ \ \ \ \ \ \ \ \ \ Junhao Shen\thanks{%
The research of second author is partially sponsored by an NSF grant.}}
\maketitle

\begin{abstract}
In this paper, we concentrate on the MF property of reduced free products of
unital C*-algebras with amalgamation over finite dimensional C*-algebras.
More specifically, we give a necessary and sufficient condition for a
reduced free product of two UHF algebras amalgamated over a
finite-dimensional C*-algebra with respect to trace preserving conditional
expectations to be MF.
\end{abstract}

\section{Introduction}

\bigskip The concept of MF algebras was first introduced by Blackadar and
Kirchberg in \cite{[BK]}. If a separable C*- algebra $\mathcal{A}$ can be
embedded into $\tprod\limits_{k}\mathcal{M}_{n_{k}}\left( \mathbb{C}\right)
/\sum_{k}\mathcal{M}_{n_{k}}\left( \mathbb{C}\right) $ for a sequence of
positive integers $\left\{ n_{k}\right\} _{k=1}^{\mathcal{1}},$ then $%
\mathcal{A}$ is an MF algebra. Many properties of MF algebras were discussed
in \cite{[BK]}. For example, it was shown there that an inductive limit of
MF algebras is an MF algebra and every subalgebra of an MF algebra is an MF
algebra. This class of C*-algebras is of interest for many reasons. For
example, it plays an important role in the classification of C*-algebras and
it is connected to the question whether the Ext semigroup of a unital
separable C*-algebra is a group. This notion is also closely connected to
Voiculescu's topological free entropy dimension for a family of self-adjoint
elements $x_{1},\cdots ,x_{n}$ in a unital C*-algebra $\mathcal{A}$ \cite%
{DV5}$.$

In \cite{DV6}, Voiculescu introduced a notion of reduced amalgamated free
product of C*-algebras in the framework of his free probability theory. The
free product construction is everywhere in the C*-literature these days, and
plays an important role in the recent study of C*-algebras \cite{IN}, \cite%
{DYKEMA}, \cite{DYKEMA1}. In \cite{HT}, Haagerup and Thorbj$\phi $rnsen
solved a long standing open problem by showing that Ext$\left( C_{r}^{\ast
}\left( F_{2}\right) \right) $ is not a group, where $C_{r}^{\ast }\left(
F_{2}\right) $ is the reduced C*-algebra of the free group $F_{2}.$ One
remarkable ingredient in their proof is their work on proving that $%
C_{r}^{\ast }\left( F_{2}\right) $ is an MF algebra. From a quick fact that $%
\left( C_{r}^{\ast }\left( F_{2}\right) ,\tau _{F_{2}}\right) =\left(
C_{r}^{\ast }\left( \mathbb{Z}\right) ,\tau _{\mathbb{Z}}\right) \ast
_{red}\left( C_{r}^{\ast }\left( \mathbb{Z}\right) ,\tau _{\mathbb{Z}%
}\right) $ where $\tau _{F_{2}},$ $\tau _{\mathbb{Z}}$ are canonical tracial
states on $C_{r}^{\ast }\left( F_{2}\right) $ and $C_{r}^{\ast }\left( 
\mathbb{Z}\right) $ respectively, one should naturally ask when a reduced
free product of two unital MF algebras is MF again. Recently, it was shown
in \cite{DS2} that a reduced free product of finite-dimensional C*-algebras
with respect to faithful tracial states is MF. More generally, the authors
proved that a reduced free product of unital separable AH algebras with
respect to faithful tracial states is an MF algebra.

In this paper, we discuss the MF property of reduced amalgamated free
products of unital C*-algebras. We investigate the reduced amalgamated free
products of finite-dimensional C*-algebras with respect to trace preserving
conditional expectations and obtain the following results. \medskip 

\textbf{Theorem 3.2 \ }Let $\mathcal{D\subseteq M}_{k}\left( \mathbb{C}%
\right) $ be a unital inclusion of C*-algebras where $\mathcal{D\ }$is a
finite-dimensional abelian C$^{\text{*}}$-algebra and $E:\mathcal{M}%
_{k}\left( \mathbb{C}\right) \rightarrow \mathcal{D}$ be the
trace-preserving conditional expectation from $\mathcal{M}_{k}\left( \mathbb{%
C}\right) \ $onto $\mathcal{D}$. Let $\left( \mathcal{M}_{k}\left( \mathbb{C}%
\right) ,E\right) \underset{\mathcal{D}}{\ast }\left( \mathcal{M}_{k}\left( 
\mathbb{C}\right) ,E\right) $ be the reduced amalgamated free product of $%
\mathcal{M}_{k}\left( \mathbb{C}\right) $ and $\mathcal{M}_{k}\left( \mathbb{%
C}\right) $ with respect to $E.$ Then $\left( \mathcal{M}_{k}\left( \mathbb{C%
}\right) ,E\right) \underset{\mathcal{D}}{\ast }\left( \mathcal{M}_{k}\left( 
\mathbb{C}\right) ,E\right) $ is an MF algebra.\medskip 

\textbf{Theorem 3.4 \ }Let $\mathcal{M}_{k}\left( \mathbb{C}\right)
\supseteq \mathcal{A\supseteq D\subseteq B\subseteq M}_{k}\left( \mathbb{C}%
\right) $ be unital C*-inclusions of finite-dimensional C*-algebras and $E:%
\mathcal{M}_{k}\left( \mathbb{C}\right) \rightarrow \mathcal{D}$ be the
trace preserving conditional expectation from $\mathcal{M}_{k}\left( \mathbb{%
C}\right) \ $onto $\mathcal{D}$. Then $\left( \mathcal{A},E|_{\mathcal{A}%
}\right) \ast _{\mathcal{D}}\left( \mathcal{B},E|_{\mathcal{B}}\right) $ is
MF.\medskip \medskip 

The following is the main result of the paper.\medskip 

\textbf{Theorem 3.5 \ }Suppose that $\mathcal{A}_{1}$ and $\mathcal{A}_{2}$
are two unital UHF-algebras with faithful tracial states $\mathcal{\tau }_{%
\mathcal{A}_{1}}$ and $\tau _{\mathcal{A}_{2}}$ respectively. Let $\mathcal{A%
}_{1}\mathcal{\supseteq D\subseteq A}_{2}$ be unital inclusions of
C*-algebras where $\mathcal{D}$ is a finite-dimensional C*-algebra. Assume
that $E_{\mathcal{A}_{1}}:\mathcal{A}_{1}\rightarrow \mathcal{D}$ and $E_{%
\mathcal{A}_{2}}:$ $\mathcal{A}_{2}\mathcal{\rightarrow D}$ are the trace
preserving conditional expectations from $\mathcal{A}_{1}$ and $\mathcal{A}%
_{2}$ onto $\mathcal{D\ }$respectively$.$ Then the reduced amalgamated free
product $\left( \mathcal{A}_{1},E_{\mathcal{A}_{1}}\right) \ast _{\mathcal{D}%
}\left( \mathcal{A}_{2},E_{\mathcal{A}_{2}}\right) $ is an MF algebra if and
only if $\tau _{\mathcal{A}_{1}}\left( z\right) =\tau _{\mathcal{A}%
_{2}}\left( z\right) $ for all $z\in \mathcal{D}.\bigskip $

A brief overview of this paper is as follows. In Section 2, we recall the
definitions of MF algebras and reduced amalgamated free product of unital
C*-algebras Section 3 is devoted to results on reduced amalgamated free
products of C*-algebras. We first show that a reduced free product of full
matrix algebras with amalgamation over a full matrix subalgebra is MF.
Later, we obtain some results about reduced amalgamated free products of two
finite-dimensional C*-algebras with amalgamation over finite-dimensional
C*-algebras. Finally, we give a necessary and sufficient condition for a
reduced free product of two UHF algebras with amalgamated over a common
finite-dimensional C*-subalgebra with respect to trace preserving
conditional expectations to be MF.

\section{Definitions and Preliminaries}

\subsection{Definition of MF Algebras}

Suppose $\{\mathcal{M}_{k_{n}}(\mathbb{C})\}_{n=1}^{\infty }$ is a sequence
of complex matrix algebras. We introduce the C*-direct product $%
\prod_{m=1}^{\infty }\mathcal{M}_{k_{m}}(\mathbb{C)}$ of $\{\mathcal{M}%
_{k_{n}}(\mathbb{C})\}_{n=1}^{\infty }$ as follows: 
\begin{equation*}
\prod_{n=1}^{\infty }\mathcal{M}_{k_{n}}(\mathbb{C})=\{(Y_{n})_{n=1}^{\infty
}\ |\ \forall \ n\geq 1,\ Y_{n}\in \mathcal{M}_{k_{n}}(C)\ \text{ and }\
\left\Vert (Y_{n})_{n=1}^{\infty }\right\Vert =\sup_{n\geq 1}\Vert
Y_{n}\Vert <\infty \}.
\end{equation*}%
Furthermore, we can introduce a norm-closed two sided ideal in $%
\prod_{n=1}^{\infty }\mathcal{M}_{k_{n}}(\mathbb{C})$ as follows: 
\begin{equation*}
\overset{\infty }{\underset{n=1}{\sum }}\mathcal{M}_{k_{n}}(\mathbb{C}%
)=\left\{ \left( Y_{n}\right) _{n=1}^{\infty }\in \prod_{n=1}^{\infty }%
\mathcal{M}_{k_{n}}(\mathbb{C}):\lim\limits_{n\rightarrow \infty }\left\Vert
Y_{n}\right\Vert =0\right\} .
\end{equation*}%
Let $\pi $ be the quotient map from $\prod_{n=1}^{\infty }\mathcal{M}%
_{k_{n}}(\mathbb{C})$ to $\prod_{n=1}^{\infty }\mathcal{M}_{k_{n}}(\mathbb{C}%
)/\overset{\infty }{\underset{n=1}{\sum }}\mathcal{M}_{k_{n}}(\mathbb{C})$.
Then 
\begin{equation*}
\prod_{n=1}^{\infty }\mathcal{M}_{k_{n}}(\mathbb{C})/\overset{\infty }{%
\underset{n=1}{\sum }}\mathcal{M}_{k_{n}}(\mathbb{C})
\end{equation*}%
is a unital C*-algebra. If we denote $\pi \left( \left( Y_{n}\right)
_{n=1}^{\infty }\right) $ by $\left[ \left( Y_{n}\right) _{n}\right] $, then 
\begin{equation*}
\left\Vert \left[ \left( Y_{n}\right) _{n}\right] \right\Vert =\underset{%
n\rightarrow \mathcal{1}}{\lim \sup }\left\Vert Y_{n}\right\Vert .
\end{equation*}%
Now we are ready to recall an equivalent definition of MF algebras given by
Blackadar and Kirchberg [\ref{[BK]}].

\begin{definition}
\label{MF}(Theorem 3.2.2, [\ref{[BK]}]) Let $\mathcal{A}$ be a separable
C*-algebra. If $\mathcal{A}$ can be embedded as a C*-subalgebra of $%
\prod_{n=1}^{\infty }\mathcal{M}_{k_{n}}(\mathbb{C})/\overset{\infty }{%
\underset{n=1}{\sum }}\mathcal{M}_{k_{n}}(\mathbb{C})$ for a sequence $%
\left\{ k_{n}\right\} _{n=1}^{\mathcal{1}}$ of integers$,$ then $\mathcal{A}$
has MF property.
\end{definition}

The examples of MF algebras contain all finite dimensional C*-algebras, AF
(approximately finite dimensional) algebras and quasidiagonal C*-algebras.
In [\ref{HT}], Haagerup and Thorbj$\phi $rnsen showed that $C_{r}^{\ast
}\left( F_{n}\right) $ is an MF algebra for $n\geq 2.$ For more examples of
MF algebras, we refer the reader to [\ref{[BK]}] and [\ref{DS}].

\subsection{Definition of Reduced Amalgamated Free Product of C*-algebras}

In this section, we will explain the KSGNS representation of unital
conditional expectations and recall a characterization of reduced
amalgamated free product of C*-algebras. First, we give the definition of
conditional expectation as follows:

\begin{definition}
\label{34}Let $\mathcal{D}$ be a unital C*-subalgebra of a unital C*-algebra 
$\mathcal{A}.$ A projection from $\mathcal{A}$ onto $\mathcal{D}$ is a
linear map $E:\mathcal{A\rightarrow D}$ such that $E\left( d\right) =d$ for
every $d\in \mathcal{D}.$ A conditional expectation from $\mathcal{A}$ onto $%
\mathcal{D}$ is a contractive completely positive projection $E$ from $%
\mathcal{A}$ onto $\mathcal{D}$ such that $E\left( dxd^{\prime }\right)
=dE\left( x\right) d^{\prime }$ for every $x\in \mathcal{A}$ and $%
d,d^{\prime }\in \mathcal{D}$ (i.e. $E$ is a $\mathcal{D}$-bimodule map).
\end{definition}

Now, we present a summary of the KSGNS representation of C*-algebras. Let$%
\mathcal{\ D}$ be a unital C*-algebra and let $\mathcal{A}$ be a unital
C*-algebra, which contains a copy of $\mathcal{D}$ as a unital
C*-subalgebra. We also suppose that there is a conditional expectation $E:%
\mathcal{A}\rightarrow \mathcal{D},$ satisfying 
\begin{equation*}
\forall a\in \mathcal{A\ }\text{with }a\neq 0,\ \exists x\in \mathcal{A}\ 
\text{such that }E\left( x^{\ast }a^{\ast }ax\right) \neq 0\ \ \ \ \ \ \ \ \
\ \ \ \left( 1\right) .
\end{equation*}%
The conditional expectations which satisfy $\left( 1\right) $ are called
non-degenerate. Let $\mathcal{H}=L^{2}\left( \mathcal{A},E\right) $ be the
right Hilbert $\mathcal{D}$-module obtained from $\mathcal{A}$ by separation
and completion with respect to the norm $\left\Vert a\right\Vert =\left\Vert
\left\langle a,a\right\rangle _{\mathcal{H}}\right\Vert ^{1/2},$ where $%
\left\langle a_{1},a_{2}\right\rangle _{\mathcal{H}}=E\left( a_{1}^{\ast
}a_{2}\right) .$ Then the linear space $\mathcal{L}\left( \mathcal{H}\right) 
$ of all adjointable $\mathcal{D}$-module operators on $\mathcal{H}$ is
actually a C*-algebra and we have a representation $\pi :\mathcal{A}%
\rightarrow \mathcal{L}\left( \mathcal{H}\right) $ defined by $\pi \left(
a\right) \widehat{a^{\prime }}=\widehat{aa^{\prime }},$ where by $\widehat{a}
$ we denote the element of $\mathcal{H},$ corresponding to $a\in \mathcal{A}.
$ $\pi $ is faithful if $E$ satisfies condition $\left( 1\right) .$ In this
construction, we have the specified element $\xi =\widehat{I_{A}}\in 
\mathcal{H}.$ We call the triple $\left( \pi ,\mathcal{H},\xi \right) $ the
KSGNS representation of $\left( \mathcal{A},E\right) $. If $\mathcal{D=}%
\mathbb{C},$ then the Hilbert C*--module considered above become simply
Hilbert space.

Let $I$ be an index set. Card$\left( I\right) \geq 2.$ Let$\mathcal{\ D}$ be
a unital C*-algebra, and for each $i\in I$, we have a unital C*-algebra $%
\mathcal{A}_{i}$, which contains a copy of $\mathcal{D}$ as a unital
C*-subalgebra. We also suppose that, for each $i\in I,$ $E_{i}:\mathcal{A}%
_{i}\rightarrow \mathcal{D},$ is non-degenerate. Let $\left( \pi _{i},%
\mathcal{H}_{i},\xi _{i}\right) $ denote the KSGNS representations of $%
\left( \mathcal{A}_{i},E_{i}\right) .$ Note that each $\pi _{i}$ is
faithful. The Hilbert $\mathcal{D}$-submodule $\xi _{i}\mathcal{D\subseteq H}%
_{i}$ is isomorphic to $\mathcal{D}$ and $\mathcal{H}_{i}\cong \xi
_{i}\oplus \mathcal{H}_{i}^{\circ }.$ Recall that $\mathcal{A}_{i}^{\circ }=%
\mathcal{A}_{i}\cap \ker E_{i}$ satisfies $\mathcal{DA}_{i}^{\circ }\mathcal{%
D\subseteq A}_{i}^{\circ }$ and $\mathcal{H}_{i}^{\circ }=\overline{\mathcal{%
A}_{i}^{\circ }\xi _{i}}$ is a C*-correspondence over $\mathcal{D}.$ We
define the free product Hilbert $\mathcal{D}$-module $\left( \mathcal{H},\xi
\right) =\ast \left( \mathcal{H}_{i},\xi _{i}\right) $ by 
\begin{equation*}
\mathcal{H=\xi D\oplus }\tbigoplus\limits_{n\geq
1}\tbigoplus\limits_{i_{1}\neq \cdots \neq i_{n}}\mathcal{H}_{i_{1}}^{\circ
}\otimes _{\mathcal{D}}\cdots \otimes _{D}\mathcal{H}_{i_{n}}^{\circ }.
\end{equation*}%
Here, $\xi \mathcal{D}$ is the trivial Hilbert $\mathcal{D}$-module $%
\mathcal{D}$ with $\xi =\widehat{1}$ and $\otimes _{\mathcal{D}}$ means
interior tensor product (see Lance's book \cite{CE}). For each $i\in I,$ let 
\begin{equation*}
\mathcal{H}\left( i\right) \mathcal{=\xi D\oplus }\tbigoplus\limits_{n\geq
1}\tbigoplus\limits_{i\neq i_{1}\neq \cdots \neq i_{n}}\mathcal{H}%
_{i_{1}}^{\circ }\otimes _{\mathcal{D}}\cdots \otimes _{D}\mathcal{H}%
_{i_{n}}^{\circ }
\end{equation*}%
and define an isomorphism $U_{i}\in \mathcal{L}\left( \mathcal{H}_{i}\otimes
_{\mathcal{D}}\mathcal{H}\left( i\right) ,\mathcal{H}\right) $ by

\begin{equation*}
U_{i}:%
\begin{array}{c}
\xi _{i}\mathcal{D\otimes }_{\mathcal{D}}\xi \mathcal{D} \\ 
\mathcal{H}_{i}^{\circ }\otimes _{\mathcal{D}}\xi \mathcal{D} \\ 
\xi _{i}\mathcal{D\otimes }_{\mathcal{D}}\left( \mathcal{H}_{i_{1}}^{\circ
}\otimes _{\mathcal{D}}\cdots \otimes _{\mathcal{D}}\mathcal{H}%
_{i_{n}}^{\circ }\right) \\ 
\mathcal{H}_{i}^{\circ }\otimes _{\mathcal{D}}\left( \mathcal{H}%
_{i_{1}}^{\circ }\otimes _{\mathcal{D}}\cdots \otimes _{\mathcal{D}}\mathcal{%
H}_{i_{n}}^{\circ }\right)%
\end{array}%
\overset{\cong }{\longrightarrow }%
\begin{array}{c}
\mathcal{\xi D} \\ 
\mathcal{H}_{i}^{\circ } \\ 
\mathcal{H}_{i_{1}}^{\circ }\otimes _{\mathcal{D}}\cdots \otimes _{\mathcal{D%
}}\mathcal{H}_{i_{n}}^{\circ } \\ 
\mathcal{H}_{i}^{\circ }\otimes _{\mathcal{D}}\mathcal{H}_{i_{1}}^{\circ
}\otimes _{\mathcal{D}}\cdots \otimes _{\mathcal{D}}\mathcal{H}%
_{i_{n}}^{\circ }%
\end{array}%
\end{equation*}%
We define a *-representation $\lambda _{i}:\mathcal{A}_{i}\rightarrow 
\mathcal{L}\left( \mathcal{H}\right) $ by $\lambda _{i}\left( x\right)
=U_{i}\left( \pi _{i}\otimes 1\right) \left( x\right) U_{i}^{\ast }.$
Observe that $\lambda _{i}|_{\mathcal{D}}$ is the canonical left action of $%
\mathcal{D}$ on the C*-correspondence $\mathcal{H}$ over $\mathcal{D},$ and
hence $\lambda _{i}|_{\mathcal{D}}=\lambda _{j}|_{\mathcal{D}}$ for every $i$
and $j.$

\begin{definition}
The reduced amalgamated free product $\left( \mathcal{A},E\right) =\ast _{%
\mathcal{D}}\left( \mathcal{A}_{i},E_{i}\right) $ is the C*-subalgebra of $%
\mathcal{L}\left( \mathcal{H}\right) $ generated by $\cup _{i\in I}\lambda
_{i}\left( \mathcal{A}_{i}\right) .$ \ 
\end{definition}

The following lemma characterizes the reduced amalgamated free product of
C*-algebras. It is one of the key ingredients in the proof of our main
result.

\begin{theorem}
(Theorem 7.2 \cite{BO}) \label{37}Let $I\in \mathcal{D\subseteq A}_{i}$ be
unital C*-algebras with non-degenerate conditional expectations $E_{i}$ from 
$\mathcal{A}_{i}$ onto $\mathcal{D},$ and let $\left( \mathcal{A},E\right)
=\ast _{\mathcal{D}}\left( \mathcal{A}_{i},E_{i}\right) .$

\begin{enumerate}
\item There is an inclusion $I\in \mathcal{D\subseteq A}$ and a
non-degenerate conditional expectation $E:\mathcal{A\rightarrow D}.$

\item We have inclusions $\mathcal{D\subseteq A}_{i}\subseteq \mathcal{A}$
which are compatible on $\mathcal{D}$, and $\mathcal{A}$ is generated by $%
\cup \mathcal{A}_{i}$ as a C*-algebra.

\item We have $E|_{\mathcal{A}_{i}}=E_{i}$ for every i, and the
C*-subalgebras $\mathcal{A}_{i}$ are free over $\mathcal{D}$ in $\left( 
\mathcal{A},E\right) .$ Namely, $E\left( a_{1}\cdots a_{n}\right) =0$ for
every $a_{j}\in $ $\mathcal{A}_{i_{j}}^{\circ }$ with $i_{1}\neq \cdots \neq
i_{n}.$
\end{enumerate}

Moreover, the above conditions uniquely characterize the reduced amalgamated
free product $\left( \mathcal{A},E\right) .$
\end{theorem}

\section{Reduced Amalgamated Free Products Of Matrix Algebras}

\subsection{$\mathcal{D}$ is a full matrix algebra}

In this subsection, we consider a reduced free product of two matrix
algebras $\mathcal{A}$ and $\mathcal{B}$ with amalgamation over a
C*-subalgebra $\mathcal{D}$ when $\mathcal{D}$ is a full matrix algebra. The
following lemma is well-known.

\begin{lemma}
\label{38.1}(Theorem 2.1, [\ref{DY}]) Let $\mathcal{A\supseteq D\subseteq B}$
be unital inclusions of unital C*-algebras. Assume $E_{1}:\mathcal{%
A\rightarrow D}$ and $E_{2}:\mathcal{B\rightarrow D}$ are non-degenerate
conditional expectations. Let $\mathcal{M}=\left( \mathcal{A},E_{1}\right)
\ast _{\mathcal{D}}\left( \mathcal{B},E_{2}\right) $ be the reduced
amalgamated free product of $\mathcal{A}$ and $\mathcal{B}$ with respect to $%
E_{1}$ and $E_{2}$. Suppose that there is a projection $p\in \mathcal{D}$
and there are partial isometries $v_{1},\cdots ,v_{n}\in \mathcal{D}$ such
that $v_{i}^{\ast }v_{i}\leq p$ and $\sum_{i=1}^{n}v_{i}v_{i}^{\ast }=1-p.$
Then%
\begin{equation*}
p\mathcal{M}p=\left( p\mathcal{A}p,E_{1}|_{p\mathcal{A}p}\right) \ast _{p%
\mathcal{D}p}\left( p\mathcal{B}p,E_{2}|_{p\mathcal{B}p}\right) .
\end{equation*}
\end{lemma}

By emulating the argument in the proof of Lemma 2.1 in [\ref{[BD]}], we
obtain the following result.

\begin{lemma}
\label{13.3}Let $\mathcal{A}$ be a unital C*-algebra and suppose there is a
projection $p\in \mathcal{A}$ and there are partial isometries $v_{1},\cdots
,v_{n}\in \mathcal{A}$ such that $v_{i}^{\ast }v_{i}\leq p$ and $\overset{n}{%
\underset{i=1}{\sum }}v_{i}v_{i}^{\ast }=1-p.$ Then $\mathcal{A}$ is MF if
and only if $p\mathcal{A}p$ is MF.
\end{lemma}

The following lemma can be found in \cite{DS2}.

\begin{lemma}
\label{40}(Theorem 3.4.1, \cite{DS2}) Suppose that $\left\{ \mathcal{A}%
_{i}:i=1,\cdots ,n\right\} $ is a family of unital separable $AH$
(approximately homogeneous) algebras with faithful tracial state $\tau _{i}$
for $i\in \left\{ 1,\cdots ,n\right\} .$ Then 
\begin{equation*}
\left( \mathcal{A}_{1},\tau _{1}\right) \ast _{red}\cdots \ast _{red}\left( 
\mathcal{A}_{n},\tau _{n}\right)
\end{equation*}%
is an MF algebra.
\end{lemma}

Now, we are able to obtain our first result on the reduced amalgamated free
products of matrix algebras.

\begin{theorem}
\label{41}Suppose that $k,n,m\in \mathbb{N}.$ If $\mathcal{M}_{k}\left( 
\mathbb{C}\right) $ can be unitally embedded into $\mathcal{M}_{n}\left( 
\mathbb{C}\right) $ and $\mathcal{M}_{m}\left( \mathbb{C}\right) $ as a
unital C*-subalgebra respectively, then for the trace preserving conditional
expectations $E_{1}$ and $E_{2}$ from $\mathcal{M}_{n}\left( \mathbb{C}%
\right) $ and $\mathcal{M}_{m}\left( \mathbb{C}\right) $ onto $\mathcal{M}%
_{k}\left( \mathbb{C}\right) $, 
\begin{equation*}
\left( \mathcal{M}_{n}\left( \mathbb{C}\right) ,E_{1}\right) \ast _{\mathcal{%
M}_{k}\left( \mathbb{C}\right) }\left( \mathcal{M}_{m}\left( \mathbb{C}%
\right) ,E_{2}\right) 
\end{equation*}%
is an MF algebra.
\end{theorem}

\begin{proof}
Since $\mathcal{M}_{k}\left( \mathbb{C}\right) $ can be unitally embedded
into $\mathcal{M}_{n}\left( \mathbb{C}\right) $ and $\mathcal{M}_{m}\left( 
\mathbb{C}\right) \ $respectively$,$ by Lemma 6.6.3 in \cite{KR}, we have
that 
\begin{equation*}
\mathcal{M}_{n}\left( \mathbb{C}\right) \cong \mathcal{M}_{k}\left( \mathbb{C%
}\right) \otimes \mathcal{M}_{\frac{n}{k}}\left( \mathbb{C}\right) 
\end{equation*}%
and 
\begin{equation*}
\mathcal{M}_{m}\left( \mathbb{C}\right) \cong \mathcal{M}_{k}\left( \mathbb{C%
}\right) \otimes \mathcal{M}_{\frac{m}{k}}\left( \mathbb{C}\right) .
\end{equation*}%
This implies that, for a one-dimensional projection $p\in \mathcal{M}%
_{k}\left( \mathbb{C}\right) ,$ we have 
\begin{equation*}
p\mathcal{M}_{n}\left( \mathbb{C}\right) p\cong \mathcal{M}_{\frac{n}{k}%
}\left( \mathbb{C}\right) \ \text{and\ }p\mathcal{M}_{m}\left( \mathbb{C}%
\right) p\cong \mathcal{M}_{\frac{m}{k}}\left( \mathbb{C}\right) .
\end{equation*}%
Suppose $\tau _{1}$ and $\tau _{2}$ are the faithful tracial states on $%
\mathcal{M}_{\frac{n}{k}}\left( \mathbb{C}\right) $ and $\mathcal{M}_{\frac{n%
}{k}}\left( \mathbb{C}\right) \ $respectively$.$ Let $E_{1}=id_{k}\otimes
\tau _{1}$ and $E_{2}=id_{k}\otimes \tau _{2}$ where $id_{k}$ is the
identity map on $\mathcal{M}_{k}\left( \mathbb{C}\right) .$ Then $E_{1}$ and 
$E_{2}$ are conditional expectations from $\mathcal{M}_{k}\left( \mathbb{C}%
\right) \otimes \mathcal{M}_{\frac{n}{k}}\left( \mathbb{C}\right) $ and $%
\mathcal{M}_{k}\left( \mathbb{C}\right) \otimes \mathcal{M}_{\frac{m}{k}%
}\left( \mathbb{C}\right) $ onto $\mathcal{M}_{k}\left( \mathbb{C}\right) $
respectively. It is easy to check that $E_{1}$ and $E_{2}$ are trace
preserving. Let 
\begin{equation*}
\mathcal{M=}\left( \mathcal{M}_{n}\left( \mathbb{C}\right) ,E_{1}\right)
\ast _{\mathcal{M}_{k}\left( \mathbb{C}\right) }\left( \mathcal{M}_{m}\left( 
\mathbb{C}\right) ,E_{2}\right) .
\end{equation*}%
Since there are partial isometries $v_{1},\cdots ,v_{k-1}\in \mathcal{M}%
_{k}\left( \mathbb{C}\right) $ such that $v_{i}^{\ast }v_{i}=p$ for each $%
i\in \left\{ 1,\cdots ,k-1\right\} $ and $\sum_{i=1}^{k-1}v_{i}v_{i}^{\ast
}=I_{\mathcal{M}_{k}\left( \mathbb{C}\right) }-p$, by Lemma \ref{38.1}, we
have 
\begin{eqnarray}
p\mathcal{M}p &=&\left( p\mathcal{M}_{n}\left( \mathbb{C}\right) p,E_{1}|_{p%
\mathcal{M}_{n}\left( \mathbb{C}\right) p}\right) \ast _{p\mathcal{M}%
_{k}\left( \mathbb{C}\right) p}\left( p\mathcal{M}_{m}\left( \mathbb{C}%
\right) p,E_{2}|_{p\mathcal{M}_{m}\left( \mathbb{C}\right) p}\right)   \notag
\\
&\cong &\left( \mathcal{M}_{\frac{n}{k}}\left( \mathbb{C}\right) ,\tau
_{1}\right) \ast _{red}\left( \mathcal{M}_{\frac{m}{k}}\left( \mathbb{C}%
\right) ,\tau _{2}\right)   \TCItag{3.1}  \label{eq3.1}
\end{eqnarray}%
Combining (3.1) with Lemma \ref{40} and Lemma \ref{13.3}, we conclude that 
\begin{equation*}
\left( \mathcal{M}_{n}\left( \mathbb{C}\right) ,E_{1}\right) \ast _{\mathcal{%
M}_{k}\left( \mathbb{C}\right) }\left( \mathcal{M}_{m}\left( \mathbb{C}%
\right) ,E_{2}\right) 
\end{equation*}%
is an MF algebra as desired.
\end{proof}

\subsection{\protect\bigskip $\mathcal{D}$ is a finite-dimensional abelian
C*-algebra }

More generally, we consider the case when the overlap C*-algebra $\mathcal{D}
$ is a finite-dimensional C*-algebra. In this subsection we consider a
reduced free product of $\mathcal{M}_{k}\left( \mathbb{C}\right) $ and $%
\mathcal{M}_{k}\left( \mathbb{C}\right) $ with amalgamation over a
finite-dimensional abelian C*-subalgebra $\mathcal{D}$ with respect to a
trace preserving conditional expectation. The following lemma can be found
in \cite{[BK]}.

\begin{lemma}
\label{40.5}(Proposition 3.3.6, \cite{[BK]}) If $\mathcal{A}$ and $\mathcal{B%
}$ are MF algebras and one of them is nuclear, then $\mathcal{A\otimes B}$
is MF.
\end{lemma}

It is well-known that the full matrix algebra $\mathcal{M}_{n}\left( \mathbb{%
C}\right) $ and the reduced C*-algebra $C_{r}^{\ast }\left( F_{n}\right) $
of the free group $F_{n}$ have canonical tracial states, without loss of
generality, we denote the reduced free product of $\mathcal{M}_{n}\left( 
\mathbb{C}\right) $ and $C_{r}^{\ast }\left( F_{n}\right) $ with respect to
their canonical tracial states by $\mathcal{M}_{n}\left( \mathbb{C}\right)
\ast _{red}C_{r}^{\ast }\left( F_{n}\right) .$ Similarly, the reduced free
product of $\mathcal{M}_{n}\left( \mathbb{C}\right) $ and $\mathcal{M}%
_{m}\left( \mathbb{C}\right) $ with respect to their tracial states is
denoted by $\mathcal{M}_{n}\left( \mathbb{C}\right) \ast _{red}\mathcal{M}%
_{m}\left( \mathbb{C}\right) .$

For the purpose of simplicity, we will first assume that $\mathcal{D}$ is a
two dimensional abelian C*-subalgebra of $\mathcal{M}_{k}\left( \mathbb{C}%
\right) $. Let $\mathcal{D}$ $\subseteq \mathcal{M}_{k}\left( \mathbb{C}%
\right) $ be a two-dimensional abelian C*-subalgebra of $\mathcal{M}%
_{k}\left( \mathbb{C}\right) ,$ i.e., there are non-zero orthogonal
projections $P_{1},$ $P_{2}$ in $\mathcal{M}_{k}\left( \mathbb{C}\right) $
such that 
\begin{equation*}
\mathcal{D=}\mathbb{C}P_{1}\oplus \mathbb{C}P_{2}.
\end{equation*}%
Let $E:\mathcal{M}_{k}\left( \mathbb{C}\right) \rightarrow \mathcal{D}$ be a
trace preserving conditional expectation from $\mathcal{M}_{k}\left( \mathbb{%
C}\right) $ onto $\mathcal{D}$ and $\left( \mathcal{M}_{k}\left( \mathbb{C}%
\right) ,E\right) \ast _{\mathcal{D}}\left( \mathcal{M}_{k}\left( \mathbb{C}%
\right) ,E\right) $ be the reduced free product of $\left( \mathcal{M}%
_{k}\left( \mathbb{C}\right) ,E\right) $ and $\left( \mathcal{M}_{k}\left( 
\mathbb{C}\right) ,E\right) $ with amalgamation over $\mathcal{D}$ with
respect to the conditional expectation $E.$ By Corollary II.6.10.8 \cite{BK}%
, $E$ is the unique trace preserving conditional expectation from $\mathcal{M%
}_{k}\left( \mathbb{C}\right) $ onto $\mathcal{D}.$ Assume that $k_{1}=$rank$%
\left( P_{1}\right) $ and $k_{2}=$rank$\left( P_{2}\right) $. Consider
C*-algebras $C_{r}^{\ast }\left( F_{2}\right) ,\ $%
\begin{equation*}
\left( \mathcal{M}_{k_{1}}\left( \mathbb{C}\right) \ast _{red}\mathcal{M}%
_{k_{2}}\left( \mathbb{C}\right) \right) \otimes \mathcal{M}_{k}\left( 
\mathbb{C}\right) 
\end{equation*}%
and 
\begin{equation*}
\left( C_{r}^{\ast }\left( F_{2}\right) \ast _{red}\mathcal{M}_{k_{1}}\left( 
\mathbb{C}\right) \ast _{red}\mathcal{M}_{k_{2}}\left( \mathbb{C}\right)
\right) \otimes \mathcal{M}_{k}\left( \mathbb{C}\right) .
\end{equation*}%
Then, by the definition of reduced free product of C*-algebras, $C_{r}^{\ast
}\left( F_{2}\right) \otimes \mathcal{M}_{k}\left( \mathbb{C}\right) \ $and $%
\left( \mathcal{M}_{k_{1}}\left( \mathbb{C}\right) \ast _{red}\mathcal{M}%
_{k_{2}}\left( \mathbb{C}\right) \right) \otimes \mathcal{M}_{k}\left( 
\mathbb{C}\right) $ can be treated as unital C*-subalgebras of 
\begin{equation*}
\left( C_{r}^{\ast }\left( F_{2}\right) \ast _{red}\mathcal{M}_{k_{1}}\left( 
\mathbb{C}\right) \ast _{red}\mathcal{M}_{k_{2}}\left( \mathbb{C}\right)
\right) \otimes \mathcal{M}_{k}\left( \mathbb{C}\right) .
\end{equation*}%
Let $U_{1}$, $U_{2}$ are canonical generators of $C_{r}^{\ast }\left(
F_{2}\right) $. Then it is easy to see that the unitary 
\begin{equation*}
U=U_{1}\otimes P_{1}+U_{2}\otimes P_{2}
\end{equation*}%
is in 
\begin{equation*}
C_{r}^{\ast }\left( F_{2}\right) \otimes \mathcal{M}_{k}\left( \mathbb{C}%
\right) \subseteq \left( C_{r}^{\ast }\left( F_{2}\right) \ast _{red}%
\mathcal{M}_{k_{1}}\left( \mathbb{C}\right) \ast _{red}\mathcal{M}%
_{k_{2}}\left( \mathbb{C}\right) \right) \otimes \mathcal{M}_{k}\left( 
\mathbb{C}\right) .
\end{equation*}

From above notations, we can get the following technical result.

\begin{lemma}
\label{43}Suppose $P_{1},P_{2}\in \mathcal{M}_{k}\left( \mathbb{C}\right) $
are orthogonal projections from above. Then there is a C*-subalgebra $%
\mathcal{M}$ of 
\begin{equation*}
\left( \mathcal{M}_{k_{1}}\left( \mathbb{C}\right) \ast _{red}\mathcal{M}%
_{k_{2}}\left( \mathbb{C}\right) \right) \otimes \mathcal{M}_{k}\left( 
\mathbb{C}\right) \subseteq \left( C_{r}^{\ast }\left( F_{2}\right) \ast
_{red}\mathcal{M}_{k_{1}}\left( \mathbb{C}\right) \ast _{red}\mathcal{M}%
_{k_{2}}\left( \mathbb{C}\right) \right) \otimes \mathcal{M}_{k}\left( 
\mathbb{C}\right) 
\end{equation*}%
such that:

\begin{enumerate}
\item $\mathcal{M}$ $\cong \mathcal{M}_{k}\left( \mathbb{C}\right) ,$

\item $\mathcal{M}_{k_{1}}\left( \mathbb{C}\right) \otimes P_{1}+\mathcal{M}%
_{k_{2}}\left( \mathbb{C}\right) \otimes P_{2}\subseteq \mathcal{M}.$
\end{enumerate}
\end{lemma}

\begin{proof}
Let $\tau _{1}$ be the faithful tracial state on $\mathcal{M}_{k_{1}}\left( 
\mathbb{C}\right) \ast _{red}\mathcal{M}_{k_{2}}\left( \mathbb{C}\right) $
obtained from faithful tracial states on $\mathcal{M}_{k_{1}}\left( \mathbb{C%
}\right) $ and $\mathcal{M}_{k_{2}}\left( \mathbb{C}\right) ,$ and $\tau _{2}
$ be the unique faithful tracial state on $\mathcal{M}_{k}\left( \mathbb{C}%
\right) $. Then the state $\tau =\tau _{1}\otimes \tau _{2}$ is a faithful
tracial state on $\left( \mathcal{M}_{k_{1}}\left( \mathbb{C}\right) \ast
_{red}\mathcal{M}_{k_{2}}\left( \mathbb{C}\right) \right) \otimes \mathcal{M}%
_{k}\left( \mathbb{C}\right) .$ Let $\left\{ e_{ij}\right\} _{i,j=1}^{k_{1}},
$ $\left\{ g_{ij}\right\} _{i,j=1}^{k_{2}}$ and $\left\{ f_{ij}\right\}
_{i,j=1}^{k}$ be systems of matrix units of $\mathcal{M}_{k_{1}}\left( 
\mathbb{C}\right) ,$ $\mathcal{M}_{k_{2}}\left( \mathbb{C}\right) $ and $%
\mathcal{M}_{k}\left( \mathbb{C}\right) ,$ respectively. So we have 
\begin{align*}
\tau \left( e_{11}\otimes P_{1}\right) & =\frac{1}{k_{1}}\cdot \frac{k_{1}}{k%
}=1\cdot \frac{1}{k} \\
& =\tau \left( I\otimes f_{11}\right) =\frac{1}{k_{2}}\cdot \frac{k_{2}}{k}
\\
& =\tau \left( g_{11}\otimes P_{2}\right) 
\end{align*}%
where $I$ is the unit in $\mathcal{M}_{k_{1}}\left( \mathbb{C}\right) \ast
_{red}\mathcal{M}_{k_{2}}\left( \mathbb{C}\right) .$ Note that $\mathcal{M}%
_{k_{1}}\left( \mathbb{C}\right) \otimes \mathcal{M}_{k}\left( \mathbb{C}%
\right) $ is *-isomorphic to a full matrix algebra. Thus there is a partial
isometry $W_{1}$ in $\mathcal{M}_{k_{1}}\left( \mathbb{C}\right) \otimes 
\mathcal{M}_{k}\left( \mathbb{C}\right) $ such that 
\begin{equation*}
W_{1}W_{1}^{\ast }=e_{11}\otimes P_{1};\ W_{1}^{\ast }W_{1}=I\otimes f_{11}.
\end{equation*}%
\newline
Similarly, there is a partial isometry $W_{2}\in \mathcal{M}_{k_{2}}\left( 
\mathbb{C}\right) \otimes \mathcal{M}_{k}\left( \mathbb{C}\right) $ such
that 
\begin{equation*}
W_{2}W_{2}^{\ast }=g_{11}\otimes P_{2};\ W_{1}^{\ast }W_{1}=I\otimes f_{11}.
\end{equation*}%
Let $V_{1}=W_{2}^{\ast }W_{1}.$ Then $V_{1}$ is a partial isometry in $%
\left( \mathcal{M}_{k_{1}}\left( \mathbb{C}\right) \ast _{red}\mathcal{M}%
_{k_{2}}\left( \mathbb{C}\right) \right) \otimes \mathcal{M}_{k}$ satisfying
\begin{equation*}
V_{1}^{\ast }V_{1}=e_{11}\otimes P_{1}\ \text{and}\ V_{1}V_{1}^{\ast
}=g_{11}\otimes P_{2}.
\end{equation*}%
Therefore the C*-subalgebra $\mathcal{M}$ in 
\begin{equation*}
\left( \mathcal{M}_{k_{1}}\left( \mathbb{C}\right) \ast _{red}\mathcal{M}%
_{k_{2}}\left( \mathbb{C}\right) \right) \otimes \mathcal{M}_{k}\left( 
\mathbb{C}\right) 
\end{equation*}%
generated by $\left\{ e_{ij}\right\} _{i,j=1}^{k_{1}}\otimes P_{1}$, $%
\left\{ g_{ij}\right\} _{i,j=1}^{k_{2}}\otimes P_{2}$ and $V_{1}$ is
*-isomorphic to the full matrix algebra $\mathcal{M}_{k}\left( \mathbb{C}%
\right) $. It is easy to see that 
\begin{equation*}
\mathcal{M}_{k_{1}}\left( \mathbb{C}\right) \otimes P_{1}+\mathcal{M}%
_{k_{2}}\left( \mathbb{C}\right) \otimes P_{2}\subseteq \mathcal{M}.
\end{equation*}
\end{proof}

\bigskip

Let $C^{\ast }\left( \mathcal{M},U\mathcal{M}U^{\ast }\right) $ be the
C*-algebra generated by $\mathcal{M}$ and $U\mathcal{M}U^{\ast }$ in 
\begin{equation*}
\left( C_{r}^{\ast }\left( F_{2}\right) \ast _{red}\mathcal{M}_{k_{1}}\left( 
\mathbb{C}\right) \ast _{red}\mathcal{M}_{k_{2}}\left( \mathbb{C}\right)
\right) \otimes \mathcal{M}_{k}\left( \mathbb{C}\right) .
\end{equation*}%
Since 
\begin{equation*}
C_{r}^{\ast }\left( F_{2}\right) \ast _{red}\mathcal{M}_{k_{1}}\left( 
\mathbb{C}\right) \ast _{red}\mathcal{M}_{k_{2}}\left( \mathbb{C}\right)
\cong C_{r}^{\ast }\left( \mathbb{Z}\right) \ast _{red}C_{r}^{\ast }\left( 
\mathbb{Z}\right) \ast _{red}\mathcal{M}_{k_{1}}\left( \mathbb{C}\right)
\ast _{red}\mathcal{M}_{k_{2}}\left( \mathbb{C}\right) ,
\end{equation*}%
we know that%
\begin{equation*}
\left( C_{r}^{\ast }\left( F_{2}\right) \ast _{red}\mathcal{M}_{k_{1}}\left( 
\mathbb{C}\right) \ast _{red}\mathcal{M}_{k_{2}}\left( \mathbb{C}\right)
\right) \otimes \mathcal{M}_{k}\left( \mathbb{C}\right)
\end{equation*}%
is an MF algebra by Lemma \ref{40.5} and \ref{40}. This implies that $%
C^{\ast }\left( \mathcal{M},U\mathcal{M}U^{\ast }\right) $ is an MF algebra.
Notice that 
\begin{equation*}
U\left( I\otimes P_{1}\right) U^{\ast }=U_{1}IU_{1}^{\ast }\otimes
P_{1}=I\otimes P_{1}
\end{equation*}%
and 
\begin{equation*}
U\left( I\otimes P_{2}\right) U^{\ast }=U_{2}IU_{2}^{\ast }\otimes
P_{2}=I\otimes P_{2},
\end{equation*}%
this implies that 
\begin{equation*}
U\left[ \mathbb{C}\left( I\otimes P_{1}\right) \oplus \mathbb{C}\left(
I\otimes P_{2}\right) \right] U^{\ast }=\mathbb{C}\left( I\otimes
P_{1}\right) \oplus \mathbb{C}\left( I\otimes P_{2}\right) .
\end{equation*}%
Therefore $\mathbb{C}\left( I\otimes P_{1}\right) \oplus \mathbb{C}\left(
I\otimes P_{2}\right) \subseteq U\mathcal{M}U^{\ast }.$

Now, we consider the reduced amalgamated free product 
\begin{equation*}
\left( \mathcal{M}_{k}\left( \mathbb{C}\right) ,E\right) \underset{\mathcal{D%
}}{\ast }\left( \mathcal{M}_{k}\left( \mathbb{C}\right) ,E\right) 
\end{equation*}%
of $\left( \mathcal{M}_{k}\left( \mathbb{C}\right) ,E\right) $ and $\left( 
\mathcal{M}_{k}\left( \mathbb{C}\right) ,E\right) $ when $\mathcal{D}$ is a
two-dimensional abelian C*-algebra.

\begin{proposition}
\label{47.5}Let $\mathcal{D\subseteq M}_{k}\left( \mathbb{C}\right) $ be a
unital inclusion of C*-algebras where $\mathcal{D}$ is a two-dimensional
abelian C*-algebra. Suppose $E:\mathcal{M}_{k}\left( \mathbb{C}\right)
\rightarrow \mathcal{D}$ is the trace preserving conditional expectation
from $\mathcal{M}_{k}\left( \mathbb{C}\right) $ onto $\mathcal{D}.$ Let $U$
and $\mathcal{M}$ be as above, then 
\begin{equation*}
C^{\ast }\left( \mathcal{M},U\mathcal{M}U^{\ast }\right) \cong \left( 
\mathcal{M}_{k}\left( \mathbb{C}\right) ,E\right) \underset{\mathcal{D}}{%
\ast }\left( \mathcal{M}_{k}\left( \mathbb{C}\right) ,E\right) .
\end{equation*}%
Therefore $\left( \mathcal{M}_{k}\left( \mathbb{C}\right) ,E\right) \underset%
{\mathcal{D}}{\ast }\left( \mathcal{M}_{k}\left( \mathbb{C}\right) ,E\right) 
$ is an MF algebra.
\end{proposition}

\begin{proof}
Let $\tau _{1}$ be the faithful tracial state on $C_{r}^{\ast }\left(
F_{2}\right) \ast _{red}\mathcal{M}_{k_{1}}\left( \mathbb{C}\right) \ast
_{red}\mathcal{M}_{k_{2}}\left( \mathbb{C}\right) $ obtained from tracial
states on $C_{r}^{\ast }\left( F_{2}\right) ,$ $\mathcal{M}_{k_{1}}\left( 
\mathbb{C}\right) $ and $\mathcal{M}_{k_{2}}\left( \mathbb{C}\right) ,$ and $%
\tau _{2}$ be the faithful tracial state on $\mathcal{M}_{k}\left( \mathbb{C}%
\right) $. Then the state $\tau =\tau _{1}\otimes \tau _{2}$ is a faithful
tracial state on 
\begin{equation*}
\left( C_{r}^{\ast }\left( F_{2}\right) \ast _{red}\mathcal{M}_{k_{1}}\left( 
\mathbb{C}\right) \ast _{red}\mathcal{M}_{k_{2}}\left( \mathbb{C}\right)
\right) \otimes \mathcal{M}_{k}\left( \mathbb{C}\right) .
\end{equation*}%
Let $I$ be the unit in $C_{r}^{\ast }\left( F_{2}\right) \ast _{red}\mathcal{%
M}_{k_{1}}\left( \mathbb{C}\right) \ast _{red}\mathcal{M}_{k_{2}}\left( 
\mathbb{C}\right) .$ It is clear that $\mathbb{C}\left( I\otimes
P_{1}\right) \oplus \mathbb{C}\left( I\otimes P_{2}\right) $ is *-isomorphic
to $\mathcal{D}.$ Assume $\mathcal{M}_{1}=\mathcal{M}$ and $\mathcal{M}_{2}=U%
\mathcal{M}U^{\ast }$. From the above discussion, we have that 
\begin{equation*}
\mathbb{C}\left( I\otimes P_{1}\right) \oplus \mathbb{C}\left( I\otimes
P_{2}\right) \subseteq \mathcal{M}_{1}\ \text{and\ }\mathbb{C}\left(
I\otimes P_{1}\right) \oplus \mathbb{C}\left( I\otimes P_{2}\right)
\subseteq \mathcal{M}_{2}.
\end{equation*}%
Define a linear mapping%
\begin{equation*}
\widetilde{E}:C^{\ast }\left( \mathcal{M}_{1},\mathcal{M}_{2}\right)
\rightarrow \mathbb{C}\left( I\otimes P_{1}\right) \oplus \mathbb{C}\left(
I\otimes P_{2}\right) 
\end{equation*}%
by 
\begin{equation*}
\widetilde{E}\left( A\right) =\sum_{i=1}^{2}\frac{1}{\tau \left( I\otimes
P_{i}\right) }\left( \tau \left( A\left( I\otimes P_{i}\right) \right)
I\right) \otimes P_{i}
\end{equation*}%
for $A\in C^{\ast }\left( \mathcal{M},U\mathcal{M}U^{\ast }\right) .$ By the
fact that 
\begin{eqnarray*}
&&I\otimes P_{1}\oplus I\otimes P_{2} \\
&=&\widetilde{I}\in \left( C_{r}^{\ast }\left( F_{2}\right) \ast _{red}%
\mathcal{M}_{k_{1}}\left( \mathbb{C}\right) \ast _{red}\mathcal{M}%
_{k_{2}}\left( \mathbb{C}\right) \right) \otimes \mathcal{M}_{k}\left( 
\mathbb{C}\right) ,
\end{eqnarray*}%
it is not hard to see that $\widetilde{E}$ is a trace preserving conditional
expectation. The restrictions of $\widetilde{E}$ on $\mathcal{M}_{1}$ and $%
\mathcal{M}_{2},$ i.e., $\widetilde{E}|_{\mathcal{M}_{1}}$ and $\widetilde{E}%
|_{\mathcal{M}_{2}}$ are trace preserving conditional expectations from $%
\mathcal{M}_{1}$ and $\mathcal{M}_{2}$ onto $\mathbb{C}\left( I\otimes
P_{1}\right) \oplus \mathbb{C}\left( I\otimes P_{2}\right) $ respectively.
Since $\mathcal{M}_{1}\cong \mathcal{M}_{k}\left( \mathbb{C}\right) $ and
the trace preserving conditional expectation from $\mathcal{M}_{k}\left( 
\mathbb{C}\right) $ onto $\mathcal{D\ }$is unique, for showing $C^{\ast
}\left( \mathcal{M}_{1},\mathcal{M}_{2}\right) \ $is *-isomorphic to $\left( 
\mathcal{M}_{k}\left( \mathbb{C}\right) ,E\right) \ast _{\mathcal{D}}\left( 
\mathcal{M}_{k}\left( \mathbb{C}\right) ,E\right) $, we only need to show
that 
\begin{equation*}
\widetilde{E}\left( T_{1}T_{2}\cdots T_{n}\right) =0
\end{equation*}%
if $T_{j}\in \mathcal{M}_{l_{j}},\ \widetilde{E}\left( T_{j}\right) =0$ for
any $j\in \left\{ 1,\cdots ,n\right\} $ and $l_{1}\neq l_{2}\neq \cdots \neq
l_{n},\ l_{1},\cdots ,l_{n}\in \left\{ 1,2\right\} $ by Theorem \ref{37}.
Suppose $X_{i},$ $Y_{i}\in \mathcal{M}_{1}$ for $i\in \left\{ 1,\cdots
,n\right\} $ with $\widetilde{E}\left( X_{i}\right) =\widetilde{E}\left(
Y_{i}\right) =0.$ If we can show that 
\begin{equation*}
\widetilde{E}\left( UX_{1}U^{\ast }Y_{1}UX_{2}U^{\ast }Y_{2}\cdots
UX_{n}U^{\ast }Y_{n}\right) =0
\end{equation*}%
for the reduced word $UX_{1}U^{\ast }Y_{1}UX_{2}U^{\ast }Y_{2}\cdots
UX_{n}U^{\ast }Y_{n}$ in $C^{\ast }\left( \mathcal{M}_{1},\mathcal{M}%
_{2}\right) ,$ then, for the other forms of reduced word, the desired result
will be followed by using similar arguments. Therefore we can conclude that $%
\left( \mathcal{M}_{k}\left( \mathbb{C}\right) ,E\right) \ast _{\mathcal{D}%
}\left( \mathcal{M}_{k}\left( \mathbb{C}\right) ,E\right) \cong C^{\ast
}\left( \mathcal{M},U\mathcal{M}U^{\ast }\right) .$ It follows that $\left( 
\mathcal{M}_{k}\left( \mathbb{C}\right) ,E\right) \ast _{\mathcal{D}}\left( 
\mathcal{M}_{k}\left( \mathbb{C}\right) ,E\right) $ is an MF algebra.

Suppose $X_{i}=X_{11}^{\left( i\right) }+X_{22}^{\left( i\right)
}+X_{12}^{\left( i\right) }+X_{21}^{\left( i\right) }$ with 
\begin{equation*}
X_{11}^{\left( i\right) }=x_{11}^{\left( i\right) }\otimes
P_{1}=\sum_{s_{1}=1,t=1}^{s_{1}=k_{1},t=k_{1}}x_{s_{1}t}^{\left( i\right)
}\otimes f_{s_{1}t}
\end{equation*}%
where $x_{11}^{\left( i\right) }\in \mathcal{M}_{k_{1}}\left( \mathbb{C}%
\right) $ and $x_{s_{1}s_{1}}^{\left( i\right) }=x_{11}^{\left( i\right) },$ 
$x_{s_{1}t}^{\left( i\right) }=0$ if $s_{1}\neq t$, and 
\begin{equation*}
X_{22}^{\left( i\right) }=x_{22}^{\left( i\right) }\otimes
P_{2}=\sum_{s_{2}=k_{1}+1,h=k_{1}+1}^{s_{2}=k,h=k}x_{s_{2}h}^{\left(
i\right) }\otimes f_{s_{2}h}
\end{equation*}%
where $x_{s_{2}s_{2}}^{\left( i\right) }=x_{22}^{\left( i\right) }\in 
\mathcal{M}_{k_{2}}\left( \mathbb{C}\right) $, $x_{s_{2}h}^{\left( i\right)
}=0$ if $s_{2}\neq h,$ and 
\begin{equation*}
X_{12}^{\left( i\right)
}=\sum_{s_{1}=1,s_{2}=k_{1}+1}^{s_{1}=k_{1},s_{2}=k}x_{s_{1}s_{2}}^{\left(
i\right) }\otimes f_{s_{1}s_{2}}\ \ \text{where\ }x_{s_{1}s_{2}}^{\left(
i\right) }\in \mathcal{M}_{k_{1}}\left( \mathbb{C}\right) \ast _{red}%
\mathcal{M}_{k_{2}}\left( \mathbb{C}\right)
\end{equation*}%
as well as 
\begin{equation*}
X_{21}^{\left( i\right)
}=\sum_{s_{1}=1,s_{2}=k_{1}+1}^{s_{1}=k_{1},s_{2}=k}x_{s_{2}s_{1}}^{\left(
i\right) }\otimes f_{s_{2}s_{1}}\ \text{where }x_{s_{2}s_{1}}^{\left(
i\right) }\in \mathcal{M}_{k_{1}}\left( \mathbb{C}\right) \ast _{red}%
\mathcal{M}_{k_{2}}\left( \mathbb{C}\right) .
\end{equation*}%
Similarly, suppose $Y_{i}=Y_{11}^{\left( i\right) }+Y_{22}^{\left( i\right)
}+Y_{12}^{\left( i\right) }+Y_{21}^{\left( i\right) }$ with%
\begin{equation*}
Y_{11}^{\left( i\right) }=y_{11}^{\left( i\right) }\otimes
P_{1}=\sum_{s_{1}=1,t=1}^{s_{1}=k_{1},t=k_{1}}y_{s_{1}t}^{\left( i\right)
}\otimes f_{s_{1}t}
\end{equation*}%
where $y_{s_{1}s_{1}}^{\left( i\right) }=y_{11}^{\left( i\right) }\in 
\mathcal{M}_{k_{1}}\left( \mathbb{C}\right) ,$ $y_{s_{1}t}^{\left( i\right)
}=0$ if $s_{1}\neq t$, and 
\begin{equation*}
Y_{22}^{\left( i\right) }=y_{22}^{\left( i\right) }\otimes
P_{2}=\sum_{s_{2}=k_{1}+1,h=k_{1}+1}^{s_{2}=k,h=k}y_{s_{2}h}^{\left(
i\right) }\otimes f_{s_{2}h}
\end{equation*}%
where $y_{s_{2}s_{2}}^{\left( i\right) }=y_{22}^{\left( i\right) }\in 
\mathcal{M}_{k_{2}}\left( \mathbb{C}\right) $, $y_{s_{2}h}^{\left( i\right)
}=0$ if $s_{2}\neq h,$ and 
\begin{equation*}
Y_{12}^{\left( i\right)
}=\sum_{s_{1}=1,s_{2}=k_{1}+1}^{s_{1}=k_{1},s_{2}=k}y_{s_{1}s_{2}}^{\left(
i\right) }\otimes f_{s_{1}s_{2}}\ \text{where }y_{s_{1}s_{2}}^{\left(
i\right) }\in \mathcal{M}_{k_{1}}\left( \mathbb{C}\right) \ast _{red}%
\mathcal{M}_{k_{2}}\left( \mathbb{C}\right)
\end{equation*}%
as well as 
\begin{equation*}
Y_{21}^{\left( i\right)
}=\sum_{s_{1}=1,s_{2}=k_{1}+1}^{s_{1}=k_{1},s_{2}=k}y_{s_{2}s_{1}}^{\left(
i\right) }\otimes f_{s_{2}s_{1}}\ \text{where }y_{s_{2}s_{1}}^{\left(
i\right) }\in \mathcal{M}_{k_{1}}\left( \mathbb{C}\right) \ast _{red}%
\mathcal{M}_{k_{2}}\left( \mathbb{C}\right) .
\end{equation*}%
Note that 
\begin{equation*}
\left( I\otimes P_{1}\right) X_{i}\left( I\otimes P_{1}\right)
=X_{11}^{\left( i\right) },\ \ \left( I\otimes P_{2}\right) X_{i}\left(
I\otimes P_{2}\right) =X_{22}^{\left( i\right) }.
\end{equation*}%
Since $\widetilde{E}\left( X_{i}\right) =0,$ we have 
\begin{equation*}
\tau \left( x_{11}^{\left( i\right) }\otimes P_{1}\right) =\tau \left(
x_{22}^{\left( i\right) }\otimes P_{1}\right) =0
\end{equation*}%
for each $i\in \left\{ 1,\cdots ,n\right\} .$ Similarly, 
\begin{equation*}
\tau \left( y_{11}^{\left( i\right) }\otimes P_{2}\right) =\tau \left(
y_{22}^{\left( i\right) }\otimes P_{2}\right) =0
\end{equation*}%
for each $i\in \left\{ 1,\cdots ,n\right\} .$ This implies that 
\begin{equation*}
\tau _{1}\left( x_{11}^{\left( i\right) }\right) =\tau _{1}\left(
x_{22}^{\left( i\right) }\right) =\tau _{1}\left( y_{11}^{\left( i\right)
}\right) =\tau _{1}\left( y_{22}^{\left( i\right) }\right) =0
\end{equation*}%
for each $i\in \left\{ 1,\cdots ,n\right\} .$ Therefore 
\begin{equation}
\tau _{1}\left( x_{s_{1}s_{1}}^{\left( i\right) }\right) =\tau _{1}\left(
x_{s_{2}s_{2}}^{\left( i\right) }\right) =\tau _{1}\left(
y_{s_{1}s_{1}}^{\left( i\right) }\right) =\tau _{1}\left(
y_{s_{2}s_{2}}^{\left( i\right) }\right) =0  \tag{3.2}  \label{eq3.2}
\end{equation}%
for each $s_{1}\in \left\{ 1,\cdots ,k_{1}\right\} $ and each $s_{2}\in
\left\{ k_{1}+1,\cdots ,k\right\} .$ Assume 
\begin{equation*}
UX_{1}U^{\ast }Y_{1}UX_{2}U^{\ast }Y_{2}\cdots UX_{n}U^{\ast
}Y_{n}=A_{11}+A_{22}+A_{12}+A_{21}
\end{equation*}%
where $A_{ij}\in \left( I\otimes P_{i}\right) C^{\ast }\left( \mathcal{M},U%
\mathcal{M}U^{\ast }\right) \left( I\otimes P_{j}\right) $ for $i,j\in
\left\{ 1,2\right\} $. Then 
\begin{align}
& A_{11}  \notag \\
& =\sum_{i_{1,}\cdots ,i_{n-1},j_{1,}\cdots ,j_{n}\in \left\{ 1,2\right\}
}\left( U_{1}\otimes P_{1}\right) X_{1j_{1}}^{\left( 1\right) }\left(
U_{j_{1}}^{\ast }\otimes P_{j_{1}}\right) Y_{j_{1}i_{1}}^{\left( 1\right)
}\cdots X_{i_{n-1}j_{n}}^{\left( n\right) }\left( U_{j_{n}}^{\ast }\otimes
P_{j_{n}}\right) Y_{j_{n}1}  \tag{3.3}  \label{eq6.1}
\end{align}%
and 
\begin{align}
& A_{22}  \notag \\
& =\sum_{i_{1,}\cdots ,i_{n-1},j_{1,}\cdots ,j_{n}\in \left\{ 1,2\right\}
}\left( U_{2}\otimes P_{2}\right) X_{2j_{1}}^{\left( 1\right) }\left(
U_{j_{1}}^{\ast }\otimes P_{j_{1}}\right) Y_{j_{1}i_{1}}\cdots
X_{i_{n-1}j_{n}}^{\left( n\right) }\left( U_{j_{n}}^{\ast }\otimes
P_{j_{n}}\right) Y_{j_{n}2}.  \tag{3.4}  \label{eq6.2}
\end{align}%
So, for showing $\widetilde{E}\left( UX_{1}U^{\ast }Y_{1}UX_{2}U^{\ast
}Y_{2}\cdots UX_{n}U^{\ast }Y_{n}\right) =0,$ we only need to prove the
following claims.\medskip

\textbf{Claim 3.1\ \ \ }For every $s_{1}\ $and every $m\ $in $\left\{
1,\cdots ,k_{1}\right\} $, every $s_{2},l\in \left\{ k_{1}+1,\cdots
,k\right\} $ and $i_{1,}\cdots ,i_{n-1},j_{1,}\cdots ,j_{n}\in \left\{
1,2\right\} ,$ we have%
\begin{equation*}
\tau _{1}\left( U_{1}x_{ms_{j_{1}}}^{\left( 1\right) }U_{j_{1}}^{\ast
}y_{s_{j_{1}}s_{i_{1}}}^{\left( 1\right) }\cdots
U_{i_{n-1}}x_{s_{i_{n-1}}s_{j_{n}}}^{\left( n\right) }U_{j_{n}}^{\ast
}y_{s_{j_{n}}m}^{\left( n\right) }\right) =0
\end{equation*}%
and \medskip 
\begin{equation*}
\tau _{1}\left( U_{2}x_{ls_{j_{1}}}^{\left( 1\right) }U_{j_{1}}^{\ast
}y_{s_{j_{1}}s_{i_{1}}}^{\left( 1\right) }\cdots
U_{i_{n-1}}x_{s_{i_{n-1}}s_{j_{n}}}^{\left( n\right) }U_{j_{n}}^{\ast
}y_{s_{j_{n}}l}^{\left( n\right) }\right) =0
\end{equation*}

\textbf{Proof} \ \ \ \ For word $U_{1}x_{ms_{j_{1}}}^{\left( 1\right)
}U_{j_{1}}^{\ast }y_{s_{j_{1}}s_{i_{1}}}^{\left( 1\right) }\cdots
U_{i_{n-1}}x_{s_{i_{n-1}}s_{j_{n}}}^{\left( n\right) }U_{j_{n}}^{\ast
}y_{s_{j_{n}}m}^{\left( n\right) },$ let $k$ denote the total number of
elements in 
\begin{equation*}
\left\{ x_{ms_{j_{1}}}^{\left( 1\right) },\cdots
,x_{s_{i_{n-1}}s_{j_{n}}}^{\left( n\right) },y_{s_{j_{1}}s_{i_{1}}}^{\left(
1\right) },\cdots ,y_{s_{j_{n}}s_{m}}^{\left( n\right) }\right\} \subseteq 
\mathcal{M}_{k_{1}}\ast \mathcal{M}_{k_{2}}
\end{equation*}%
which have non-zero trace value. If $k=0$, then each element in 
\begin{equation*}
\left\{ x_{ms_{j_{1}}}^{\left( 1\right) },\cdots
,x_{s_{i_{n-1}}s_{j_{n}}}^{\left( n\right) },y_{s_{j_{1}}s_{i_{1}}}^{\left(
1\right) },\cdots ,y_{s_{j_{n}}s_{m}}^{\left( n\right) }\right\} ,
\end{equation*}%
has zero trace value. Since $U_{1}$ and $U_{2}$ are two free Haar unitaries
in $C_{r}^{\ast }\left( F_{2}\right) $ and 
\begin{equation*}
U_{1}x_{ms_{j_{1}}}^{\left( 1\right) }U_{j_{1}}^{\ast
}y_{s_{j_{1}}s_{i_{1}}}^{\left( 1\right) }\cdots
U_{i_{n-1}}x_{s_{i_{n-1}}s_{j_{n}}}^{\left( n\right) }U_{j_{n}}^{\ast
}y_{s_{j_{n}}m}^{\left( n\right) }
\end{equation*}%
is a reduced word in $C_{r}^{\ast }\left( F_{2}\right) \ast _{red}\left( 
\mathcal{M}_{k_{1}}\left( \mathbb{C}\right) \ast _{red}\mathcal{M}%
_{k_{2}}\left( \mathbb{C}\right) \right) ,$ we have 
\begin{equation*}
\tau _{1}\left( U_{1}x_{ms_{j_{1}}}^{\left( 1\right) }U_{j_{1}}^{\ast
}y_{s_{j_{1}}s_{i_{1}}}^{\left( 1\right) }\cdots
U_{i_{n-1}}x_{s_{i_{n-1}}s_{j_{n}}}^{\left( n\right) }U_{j_{n}}^{\ast
}y_{s_{j_{n}}m}^{\left( n\right) }\right) =0.
\end{equation*}%
Assume that $k=1.$ Then there is an integer $t$ such that $\tau _{1}\left(
x_{s_{i_{_{t-1}}}s_{j_{t}}}^{\left( t\right) }\right) $ or $\tau _{1}\left(
y_{s_{j_{t}}s_{i_{t}}}^{\left( t\right) }\right) $ is non-zero. From
equation (3.2), we get that $i_{t-1}\neq j_{t}$ or $j_{t}\neq i_{t}.$
Without loss of generality, we may assume that $\tau _{1}\left(
y_{s_{j_{t}}s_{i_{t}}}^{\left( t\right) }\right) \neq 0$, $1=j_{t}\neq
i_{t}=2$ and let $\overset{\circ }{y_{s_{1}s_{2}}^{\left( t\right) }}%
=y_{s_{1}s_{2}}^{\left( t\right) }-\tau _{1}\left( y_{s_{1}s_{2}}^{\left(
t\right) }\right) .$ Therefore 
\begin{align*}
& U_{1}x_{ms_{j_{1}}}^{\left( 1\right) }U_{j_{1}}^{\ast
}y_{s_{j_{1}}s_{i_{1}}}^{\left( 1\right) }\cdots
U_{i_{n-1}}x_{s_{i_{n-1}}s_{j_{n}}}^{\left( n\right) }U_{j_{n}}^{\ast
}y_{s_{j_{n}}m}^{\left( n\right) } \\
& =U_{1}x_{ms_{j_{1}}}^{\left( 1\right) }U_{j_{1}}^{\ast
}y_{s_{j_{1}}s_{i_{1}}}^{\left( 1\right) }\cdots U_{1}^{\ast }\overset{\circ 
}{y_{s_{1}s_{2}}^{\left( t\right) }}U_{2}\cdots
U_{i_{n-1}}x_{s_{i_{n-1}}s_{j_{n}}}^{\left( n\right) }U_{j_{n}}^{\ast
}y_{s_{j_{n}}m}^{\left( n\right) } \\
& +\tau _{1}\left( y_{s_{1}s_{2}}^{\left( t\right) }\right)
U_{1}x_{ms_{j_{1}}}^{\left( 1\right) }U_{j_{1}}^{\ast
}y_{s_{j_{1}}s_{i_{1}}}^{\left( 1\right) }\cdots U_{1}^{\ast }U_{2}\cdots
U_{i_{n-1}}x_{s_{i_{n-1}}s_{j_{n}}}^{\left( n\right) }U_{j_{n}}^{\ast
}y_{s_{j_{n}}s_{m}}^{\left( n\right) }
\end{align*}%
It is easy to see that 
\begin{equation*}
U_{1}x_{m,s_{j_{1}}}^{\left( 1\right) }U_{j_{1}}^{\ast
}y_{s_{j_{1}}s_{i_{1}}}^{\left( 1\right) }\cdots U_{1}^{\ast }\overset{\circ 
}{y_{s_{1}s_{2}}^{\left( t\right) }}U_{2}\cdots
U_{i_{n-1}}x_{s_{i_{n-1}}s_{j_{n}}}^{\left( n\right) }U_{j_{n}}^{\ast
}y_{s_{j_{n}}m}^{\left( n\right) }
\end{equation*}%
and 
\begin{equation*}
U_{1}x_{m,s_{j_{1}}}^{\left( 1\right) }U_{j_{1}}^{\ast
}y_{s_{j_{1}}s_{i_{1}}}^{\left( 1\right) }\cdots U_{1}^{\ast }U_{2}\cdots
U_{i_{n-1}}x_{s_{i_{n-1}}s_{j_{n}}}^{\left( n\right) }U_{j_{n}}^{\ast
}y_{s_{j_{n}}s_{m}}^{\left( n\right) }
\end{equation*}%
are reduced words in $C_{r}^{\ast }\left( F_{2}\right) \ast _{red}\left( 
\mathcal{M}_{k_{1}}\left( \mathbb{C}\right) \ast _{red}\mathcal{M}%
_{k_{2}}\left( \mathbb{C}\right) \right) $. Since $k=1,$ we have the total
number of elements in each reduced word from above which belong to $\mathcal{%
M}_{k_{1}}\ast \mathcal{M}_{k_{2}}$ with non-zero trace is $0.$ It follows
that 
\begin{eqnarray*}
&&\tau _{1}\left( U_{1}x_{ms_{j_{1}}}^{\left( 1\right) }U_{j_{1}}^{\ast
}y_{s_{j_{1}}s_{i_{1}}}^{\left( 1\right) }\cdots U_{1}^{\ast }\overset{\circ 
}{y_{s_{1}s_{2}}^{\left( t\right) }}U_{2}\cdots
U_{i_{n-1}}x_{s_{i_{n-1}}s_{j_{n}}}^{\left( n\right) }U_{j_{n}}^{\ast
}y_{s_{j_{n}}m}^{\left( n\right) }\right)  \\
&=&\tau _{1}\left( U_{1}x_{ms_{j_{1}}}^{\left( 1\right) }U_{j_{1}}^{\ast
}y_{s_{j_{1}}s_{i_{1}}}^{\left( 1\right) }\cdots U_{1}^{\ast }U_{2}\cdots
U_{i_{n-1}}x_{s_{i_{n-1}}s_{j_{n}}}^{\left( n\right) }U_{j_{n}}^{\ast
}y_{s_{j_{n}}m}^{\left( n\right) }\right) =0.
\end{eqnarray*}%
By a similar argument, if $k\geq 2,$ we can decompose 
\begin{equation*}
U_{1}x_{ms_{j_{1}}}^{\left( 1\right) }U_{j_{1}}^{\ast
}y_{s_{j_{1}}s_{i_{1}}}^{\left( 1\right) }\cdots
U_{i_{n-1}}x_{s_{i_{n-1}}s_{j_{n}}}^{\left( n\right) }U_{j_{n}}^{\ast
}y_{s_{j_{n}}m}^{\left( n\right) }
\end{equation*}%
as a sum of two reduced words in 
\begin{equation*}
C_{r}^{\ast }\left( F_{2}\right) \ast _{red}\left( \mathcal{M}_{k_{1}}\left( 
\mathbb{C}\right) \ast _{red}\mathcal{M}_{k_{2}}\left( \mathbb{C}\right)
\right) 
\end{equation*}%
such that, for each reduced word, the total number of elements in $\mathcal{M%
}_{k_{1}}\ast \mathcal{M}_{k_{2}}$ with non-zero trace is $k-1.$ This
implies that we can decompose 
\begin{equation*}
U_{1}x_{ms_{j_{1}}}^{\left( 1\right) }U_{j_{1}}^{\ast
}y_{s_{j_{1}}s_{i_{1}}}^{\left( 1\right) }\cdots
U_{i_{n-1}}x_{s_{i_{n-1}}s_{j_{n}}}^{\left( n\right) }U_{j_{n}}^{\ast
}y_{s_{j_{n}}m}^{\left( n\right) }
\end{equation*}%
as a sum of finitely many reduced words satisfying the condition that the
total number of elements which belong to $\mathcal{M}_{k_{1}}\ast \mathcal{M}%
_{k_{2}}$ in the reduced word with non-zero trace is $0$. Therefore we have
that the trace value of each summand is $0.$ It implies that 
\begin{equation*}
\tau _{1}\left( U_{1}x_{ms_{j_{1}}}^{\left( 1\right) }U_{j_{1}}^{\ast
}y_{s_{j_{1}}s_{i_{1}}}^{\left( 1\right) }\cdots
U_{i_{n-1}}x_{s_{i_{n-1}}s_{j_{n}}}^{\left( n\right) }U_{j_{n}}^{\ast
}y_{s_{j_{n}}m}^{\left( n\right) }\right) =0.
\end{equation*}%
By the same argument, we conclude that 
\begin{equation*}
\tau _{1}\left( U_{2}x_{ls_{j_{1}}}^{\left( 1\right) }U_{j_{1}}^{\ast
}y_{s_{j_{1}}s_{i_{1}}}^{\left( 1\right) }\cdots
U_{i_{n-1}}x_{s_{i_{n-1}}s_{j_{n}}}^{\left( n\right) }U_{j_{n}}^{\ast
}y_{s_{j_{n}}l}^{\left( n\right) }\right) =0
\end{equation*}

So this complete the proof of Claim 3.1.

\textbf{Claim 3.2\ \ } If 
\begin{equation*}
\tau _{1}\left( U_{1}x_{ms_{j_{1}}}^{\left( 1\right) }U_{j_{1}}^{\ast
}y_{s_{j_{1}}s_{i_{1}}}^{\left( 1\right) }\cdots
U_{i_{n-1}}x_{s_{i_{n-1}}s_{j_{n}}}^{\left( n\right) }U_{j_{n}}^{\ast
}y_{s_{j_{n}}m}^{\left( n\right) }\right) =0\ 
\end{equation*}%
and 
\begin{equation*}
\tau _{1}\left( U_{2}x_{ls_{j_{1}}}^{\left( 1\right) }U_{j_{1}}^{\ast
}y_{s_{j_{1}}s_{i_{1}}}^{\left( 1\right) }\cdots
U_{i_{n-1}}x_{s_{i_{n-1}}s_{j_{n}}}^{\left( n\right) }U_{j_{n}}^{\ast
}y_{s_{j_{n}}l}^{\left( n\right) }\right) =0
\end{equation*}%
for each $s_{1},m\in \left\{ 1,\cdots ,k_{1}\right\} $, $s_{2},l\in \left\{
k_{1}+1,\cdots ,k\right\} $ and $i_{1,}\cdots ,i_{n-1},j_{1,}\cdots
,j_{n}\in \left\{ 1,2\right\} ,$ then 
\begin{equation*}
\widetilde{E}\left( UX_{1}U^{\ast }Y_{1}UX_{2}U^{\ast }Y_{2}\cdots
UX_{n}U^{\ast }Y_{n}\right) =0.\medskip 
\end{equation*}

\textbf{Proof \ \ \ \ }Suppose%
\begin{eqnarray*}
&&\left( U_{1}\otimes P_{1}\right) X_{1j_{1}}^{\left( 1\right) }\left(
U_{j_{1}}^{\ast }\otimes P_{j_{1}}\right) Y_{j_{1}i_{1}}^{\left( 1\right)
}\cdots \left( U_{i_{n-1}}\otimes P_{i_{n-1}}\right)
X_{i_{n-1}j_{n}}^{\left( n\right) }\left( U_{j_{n}}^{\ast }\otimes
P_{j_{n}}\right) Y_{j_{n}1}^{\left( n\right) } \\
&=&\sum_{s,t=1}^{k_{1}}T_{st}\otimes f_{st}.
\end{eqnarray*}%
in (3.3). Then 
\begin{equation}
T_{mm}=\sum_{s_{2}=k_{1}+1}^{k}\sum_{s_{1}=1}^{k_{1}}\sum_{j_{1},\cdots
,j_{n},i_{1},\cdots ,i_{n-1}\in \left\{ 1,2\right\}
}U_{1}x_{ms_{j_{1}}}^{\left( 1\right) }U_{j_{1}}^{\ast
}y_{s_{j_{1}}s_{i_{1}}}^{\left( 1\right) }\cdots
U_{i_{n-1}}x_{s_{i_{n-1}}s_{j_{n}}}^{\left( n\right) }U_{j_{n}}^{\ast
}y_{s_{j_{n}}m}^{\left( n\right) }\ \   \tag{3.5}  \label{eq3.5}
\end{equation}%
for$\ m\in \left\{ 1,\cdots ,k_{1}\right\} .$ Since 
\begin{align}
& \tau \left( \left( U_{1}\otimes P_{1}\right) X_{1j_{1}}^{\left( 1\right)
}\left( U_{j_{1}}^{\ast }\otimes P_{j_{1}}\right) Y_{j_{1}i_{1}}^{\left(
1\right) }\cdots \left( U_{i_{n-1}}\otimes P_{i_{n-1}}\right)
X_{i_{n-1}j_{n}}^{\left( n\right) }\left( U_{j_{n}}^{\ast }\otimes
P_{j_{n}}\right) Y_{j_{n}1}^{\left( n\right) }\right)   \notag \\
& =\sum_{m=1}^{k_{1}}\frac{1}{k}\tau _{1}\left( T_{mm}\right)   \tag{3.6}
\label{eq3.6}
\end{align}%
and 
\begin{equation*}
\tau _{1}\left( U_{1}x_{ms_{j_{1}}}^{\left( 1\right) }U_{j_{1}}^{\ast
}y_{s_{j_{1}}s_{i_{1}}}^{\left( 1\right) }\cdots
U_{i_{n-1}}x_{s_{i_{n-1}}s_{j_{n}}}^{\left( n\right) }U_{j_{n}}^{\ast
}y_{s_{j_{n}}m}^{\left( n\right) }\right) =0
\end{equation*}%
for each $s_{1},m\in \left\{ 1,\cdots ,k_{1}\right\} $, $s_{2}\in \left\{
k_{1}+1,\cdots ,k\right\} $ and $i_{1,}\cdots ,i_{n-1},j_{1,}\cdots
,j_{n}\in \left\{ 1,2\right\} ,$ from $\left( \ref{eq3.5}\right) $ and $%
\left( \ref{eq3.6}\right) $, we get 
\begin{equation*}
\tau \left( \left( U_{1}\otimes P_{1}\right) X_{1j_{1}}^{\left( 1\right)
}\left( U_{j_{1}}^{\ast }\otimes P_{j_{1}}\right) Y_{j_{1}i_{1}}^{\left(
1\right) }\cdots X_{i_{n-1}j_{n}}^{\left( n\right) }\left( U_{j_{n}}^{\ast
}\otimes P_{j_{n}}\right) Y_{j_{n}1}\right) =0
\end{equation*}%
Similarly, we can get that 
\begin{equation*}
\tau \left( \left( U_{2}\otimes P_{2}\right) X_{2j_{1}}^{\left( 1\right)
}\left( U_{j_{1}}^{\ast }\otimes P_{j_{1}}\right) Y_{j_{1}i_{1}}\cdots
X_{i_{n-1}j_{n}}^{\left( n\right) }\left( U_{j_{n}}^{\ast }\otimes
P_{j_{n}}\right) Y_{j_{n}2}\right) =0
\end{equation*}%
for $i_{1,}\cdots ,i_{n-1},j_{1,}\cdots ,j_{n}\in \left\{ 1,2\right\} $. It
follows that $\tau \left( A_{11}\right) =\tau \left( A_{22}\right) =0.$ Note
that%
\begin{equation*}
\widetilde{E}\left( UX_{1}U^{\ast }Y_{1}UX_{2}U^{\ast }Y_{2}\cdots
UX_{n}U^{\ast }Y_{n}\right) =\sum_{i=1}^{2}\frac{1}{\tau \left( I\otimes
P_{i}\right) }\left( \tau \left( A_{ii}\right) I\otimes P_{i}\right) ,
\end{equation*}%
therefore 
\begin{equation*}
\widetilde{E}\left( UX_{1}U^{\ast }Y_{1}UX_{2}U^{\ast }Y_{2}\cdots
UX_{n}U^{\ast }Y_{n}\right) =0
\end{equation*}%
in $C_{r}^{\ast }\left( F_{2}\right) \ast _{red}\left( \mathcal{M}%
_{k_{1}}\left( \mathbb{C}\right) \ast _{red}\mathcal{M}_{k_{2}}\left( 
\mathbb{C}\right) \right) $.\medskip 

Then, from the above argument, we conclude that $\left( \mathcal{M}%
_{k}\left( \mathbb{C}\right) ,E\right) \underset{\mathcal{D}}{\ast }\left( 
\mathcal{M}_{k}\left( \mathbb{C}\right) ,E\right) $ is an MF
algebra..\bigskip
\end{proof}

Although Proposition 3.1 is stated for a two-dimensional abelian
C*-subalgebra, similar result holds for each finite-dimensional abelian
C*-subalgebra by using a similar argument. We state this as a Theorem

\begin{theorem}
\label{48}Let $\mathcal{D\subseteq M}_{k}\left( \mathbb{C}\right) $ be a
unital inclusion of C*-algebras where $\mathcal{D}$ is a finite-dimensional
abelian C*-algebra. Then, for the trace preserving conditional expectation $%
E:\mathcal{M}_{k}\left( \mathbb{C}\right) \rightarrow \mathcal{D},$ the
reduced amalgamated free product $\left( \mathcal{M}_{k}\left( \mathbb{C}%
\right) ,E\right) \underset{\mathcal{D}}{\ast }\left( \mathcal{M}_{k}\left( 
\mathbb{C}\right) ,E\right) $ of $\left( \mathcal{M}_{k}\left( \mathbb{C}%
\right) ,E\right) $ and $\left( \mathcal{M}_{k}\left( \mathbb{C}\right)
,E\right) $ with respect to $E$ is an MF algebra.\medskip 
\end{theorem}

\textbf{Sketch of the proof}\ \ We might assume that there are mutually
orthogonal projections $P_{1},\cdots ,P_{n}$ in $\mathcal{M}_{k}\left( 
\mathbb{C}\right) $ such that 
\begin{equation*}
\mathcal{D\cong \ }\mathbb{C}P_{1}\oplus \mathbb{C}P_{2}\oplus \mathbb{%
\cdots \oplus C}P_{n}.
\end{equation*}%
Let $k_{i}=$rank$\left( P_{i}\right) $ for $1\leq i\leq n.$ Consider the
unital C*-algebra:%
\begin{equation*}
\left( C_{r}^{\ast }\left( F_{n}\right) \ast _{red}\mathcal{M}_{k_{1}}\left( 
\mathbb{C}\right) \ast _{red}\cdots \ast _{red}\mathcal{M}_{k_{n}}\left( 
\mathbb{C}\right) \right) \otimes \mathcal{M}_{k}\left( \mathbb{C}\right) .
\end{equation*}%
Let $U_{1},\cdots ,U_{n}$ be the canonical generators of $C_{r}^{\ast
}\left( F_{n}\right) $ and 
\begin{equation*}
U=U_{1}\otimes P_{1}+\cdots +U_{n}\otimes P_{n}
\end{equation*}%
be a unitary in 
\begin{eqnarray*}
&&C_{r}^{\ast }\left( F_{n}\right) \otimes \mathcal{M}_{k}\left( \mathbb{C}%
\right)  \\
&\subseteq &\left( C_{r}^{\ast }\left( F_{n}\right) \ast _{red}\mathcal{M}%
_{k_{1}}\left( \mathbb{C}\right) \ast _{red}\cdots \ast _{red}\mathcal{M}%
_{k_{n}}\left( \mathbb{C}\right) \right) \otimes \mathcal{M}_{k}\left( 
\mathbb{C}\right) .
\end{eqnarray*}%
Similar argument in Lemma \ref{43}, shows that there is a unital
C*-subalgebra $\mathcal{M}$ of $\left( \mathcal{M}_{k_{1}}\left( \mathbb{C}%
\right) \ast _{red}\cdots \ast _{red}\mathcal{M}_{k_{n}}\left( \mathbb{C}%
\right) \right) \otimes \mathcal{M}_{k}\left( \mathbb{C}\right) $ such that $%
\mathcal{M\cong M}_{k}\left( \mathbb{C}\right) $ and 
\begin{equation*}
\mathcal{M}_{k_{1}}\left( \mathbb{C}\right) \otimes P_{1}+\cdots \oplus 
\mathcal{M}_{k_{n}}\left( \mathbb{C}\right) \otimes P_{n}\subseteq \mathcal{M%
}
\end{equation*}%
as well as 
\begin{equation*}
\mathbb{C}\left( I\otimes P_{1}\right) \oplus \cdots \oplus \mathbb{C}\left(
I\otimes P_{n}\right) \subseteq U\mathcal{M}U^{\ast }.
\end{equation*}%
Let $C^{\ast }\left( \mathcal{M},U\mathcal{M}U^{\ast }\right) $ be a unital
C*-subalgebra generated by $\mathcal{M}$ and $U\mathcal{M}U^{\ast }$ in $%
\left( C_{r}^{\ast }\left( F_{n}\right) \ast _{red}\mathcal{M}_{k_{1}}\left( 
\mathbb{C}\right) \ast _{red}\cdots \ast _{red}\mathcal{M}_{k_{n}}\left( 
\mathbb{C}\right) \right) \otimes \mathcal{M}_{k}\left( \mathbb{C}\right) .$
Then $C^{\ast }\left( \mathcal{M},U\mathcal{M}U^{\ast }\right) $ is an MF
algebra. Similar argument as Proposition \ref{47.5} shows that 
\begin{equation*}
C^{\ast }\left( \mathcal{M},U\mathcal{M}U^{\ast }\right) \cong \left( 
\mathcal{M}_{k}\left( \mathbb{C}\right) ,E\right) \underset{\mathcal{D}}{%
\ast }\left( \mathcal{M}_{k}\left( \mathbb{C}\right) ,E\right) .
\end{equation*}%
This induces that $\left( \mathcal{M}_{k}\left( \mathbb{C}\right) ,E\right) 
\underset{\mathcal{D}}{\ast }\left( \mathcal{M}_{k}\left( \mathbb{C}\right)
,E\right) $ is also an MF algebra, which completes our proof. $\blacksquare
\bigskip $

\subsection{$\mathcal{D}$ is a finite-dimensional C*-algebra}

Now, we are ready to consider the case when $\mathcal{D}$ is a
finite-dimensional C*-algebra.

\begin{theorem}
\label{49}Let $\mathcal{D\subseteq M}_{k}\left( \mathbb{C}\right) $ be a
unital inclusion of C*-algebras where $\mathcal{D}$ is a finite-dimensional C%
$^{\text{*}}$-algebra and $E:\mathcal{M}_{k}\left( \mathbb{C}\right)
\rightarrow \mathcal{D}$ be the trace preserving conditional expectation
from $\mathcal{M}_{k}\left( \mathbb{C}\right) $ onto $\mathcal{D}$. Then $%
\left( \mathcal{M}_{k}\left( \mathbb{C}\right) ,E\right) \ast _{\mathcal{D}%
}\left( \mathcal{M}_{k}\left( \mathbb{C}\right) ,E\right) $ is MF.
\end{theorem}

\begin{proof}
\textbf{\ }Suppose $\mathcal{D=M}_{k_{1}}\left( \mathbb{C}\right) \oplus
\cdots \oplus \mathcal{M}_{k_{t}}\left( \mathbb{C}\right) \ $with $%
k_{1}+\cdots +k_{t}=k$ and $P_{1,}\cdots ,P_{t}$ are central projections of $%
\mathcal{D}$ with $\sum_{i=1}^{t}P_{i}=I\in \mathcal{M}_{k}\left( \mathbb{C}%
\right) .$ Let $p_{i}$ be a minimal projection in $P_{i}\mathcal{D}$ for
each $i\in \left\{ 1,\cdots ,t\right\} .$ Then $p_{i}\mathcal{D}p_{i}\cong 
\mathbb{C}p_{i}.$ Assume $p=\sum_{i=1}^{t}p_{i},$ it follows that $p\mathcal{%
D}p$ is a finite-dimensional abelian C*-algebra. Note that there are partial
isometries $v_{1},\cdots ,v_{n}\in \mathcal{D}$ with $v_{i}^{\ast }v_{i}\leq
p$ and $\sum_{i=1}^{n}v_{i}v_{i}^{\ast }\leq 1-p$ where $n=\max \left\{
k_{1},\cdots ,k_{t}\right\} -1.$ So, by Lemma \ref{38.1} and Lemma \ref{13.3}%
, it is sufficient to show that 
\begin{equation*}
\left( p\mathcal{M}_{k}\left( \mathbb{C}\right) p,E|_{p\mathcal{M}_{k}\left( 
\mathbb{C}\right) p}\right) \ast _{p\mathcal{D}p}\left( p\mathcal{M}%
_{k}\left( \mathbb{C}\right) p,E|_{p\mathcal{M}_{k}\left( \mathbb{C}\right)
p}\right) 
\end{equation*}%
is an MF algebra. By Theorem \ref{48}, it is easy to see that 
\begin{equation*}
\left( p\mathcal{M}_{k}\left( \mathbb{C}\right) p,E|_{p\mathcal{M}_{k}\left( 
\mathbb{C}\right) p}\right) \ast _{p\mathcal{D}p}\left( p\mathcal{M}%
_{k}\left( \mathbb{C}\right) p,E|_{p\mathcal{M}_{k}\left( \mathbb{C}\right)
p}\right) 
\end{equation*}%
is MF, so is $\left( \mathcal{M}_{k}\left( \mathbb{C}\right) ,E\right) \ast
_{\mathcal{D}}\left( \mathcal{M}_{k}\left( \mathbb{C}\right) ,E\right) $. 
\end{proof}

It is not hard to extend the previous theorem to a reduced amalgamated free
product of two finite-dimensional C*-algebras.

\begin{theorem}
\label{70}Let $\mathcal{M}_{k}\left( \mathbb{C}\right) \supseteq \mathcal{%
A\supseteq D\subseteq B\subseteq M}_{k}\left( \mathbb{C}\right) $ be unital
C*-inclusions of finite-dimensional C*-algebras and $E:\mathcal{M}_{k}\left( 
\mathbb{C}\right) \rightarrow \mathcal{D}$ be the trace preserving
conditional expectation. Then $\left( \mathcal{A},E|_{\mathcal{A}}\right)
\ast _{\mathcal{D}}\left( \mathcal{B},E|_{\mathcal{B}}\right) $ is MF.
\end{theorem}

\begin{proof}
From the definition of reduced amalgamated free product of C*-algebras, we
have the following C*-inclusion 
\begin{equation*}
\left( \mathcal{A},E|_{\mathcal{A}}\right) \ast _{\mathcal{D}}\left( 
\mathcal{B},E|_{\mathcal{B}}\right) \subseteq \left( \mathcal{M}_{k}\left( 
\mathbb{C}\right) ,E\right) \ast _{\mathcal{D}}\left( \mathcal{M}_{k}\left( 
\mathbb{C}\right) ,E\right) .
\end{equation*}%
So $\left( \mathcal{A},E|_{\mathcal{A}}\right) \ast _{\mathcal{D}}\left( 
\mathcal{B},E|_{\mathcal{B}}\right) $ is an MF algebra.
\end{proof}

Now we are ready to prove our main result in the paper. 

Let $\mathbb{C\langle }\mathbf{X}_{1},\ldots ,\mathbf{X}_{n}\mathbb{\rangle }
$ be the set of all noncommutative polynomials in the indeterminants $%
\mathbf{X}_{1},\ldots ,\mathbf{X}_{n}$. Then we let%
\begin{equation*}
\mathbb{C\langle }\mathbf{X}_{1},\mathbf{X}_{2},\cdots \mathbb{\rangle =\cup 
}_{m=1}^{\infty }\mathbb{C\langle }\mathbf{X}_{1},\mathbf{X}_{2},\cdots 
\mathbf{X}_{m}\mathbb{\rangle }.
\end{equation*}%
The following theorem give a necessary and sufficient condition such that a
reduced free product of two UHF algebras with amalgamation over a
finite-dimensional C*-algebra is MF.

\begin{theorem}
\label{70.1}Suppose that $\mathcal{A}_{1}$ and $\mathcal{A}_{2}$ are two
unital UHF-algebras with the faithful tracial states $\mathcal{\tau }_{%
\mathcal{A}_{1}}$ and $\tau _{\mathcal{A}_{2}}$ respectively. Let $\mathcal{A%
}_{1}\mathcal{\supseteq D\subseteq A}_{2}$ be unital inclusions of
C*-algebras where $\mathcal{D}$ is a finite-dimensional C*-algebra. Assume
that $E_{\mathcal{A}_{1}}:\mathcal{A}_{1}\rightarrow \mathcal{D}$ and $E_{%
\mathcal{A}_{2}}:$ $\mathcal{A}_{2}\mathcal{\rightarrow D}$ are the trace
preserving conditional expectations onto $\mathcal{D}.$ Then the reduced
amalgamated free product $\left( \mathcal{A}_{1},E_{\mathcal{A}_{1}}\right)
\ast _{\mathcal{D}}\left( \mathcal{A}_{2},E_{\mathcal{A}_{2}}\right) $ is an
MF algebra if and only if $\tau _{\mathcal{A}_{1}}\left( z\right) =\tau _{%
\mathcal{A}_{2}}\left( z\right) $ for all $z\in \mathcal{D}.$
\end{theorem}

\begin{proof}
From the fact that every MF algebra has a tracial state and $\mathcal{A}_{1},
$ $\mathcal{A}_{2}$ have unique tracial states, one direction of the proof
is obvious. Suppose $\mathcal{A}_{1}\mathcal{=}\overline{\mathcal{\cup }%
_{k=1}^{\mathcal{1}}\mathcal{M}_{n_{k}}\left( \mathbb{C}\right) }$ and $%
\mathcal{A}_{2}\mathcal{=}\overline{\mathcal{\cup }_{l=1}^{\mathcal{1}}%
\mathcal{M}_{m_{l}}\left( \mathbb{C}\right) }$ where $\left\{ \mathcal{M}%
_{n_{k}}\left( \mathbb{C}\right) \right\} _{k=1}^{\mathcal{1}}$ and $\left\{ 
\mathcal{M}_{m_{l}}\left( \mathbb{C}\right) \right\} _{l=1}^{\mathcal{1}}$
are increasing sequences of full matrix algebras. Since $\mathcal{D}$ is a
finite-dimensional C*-algebra, without loss of generality, we may assume
that there is an integer $t_{0}$ such that, for $\forall t\geq t_{0},$ we
have $\mathcal{D\subseteq M}_{n_{t}}\left( \mathbb{C}\right) $ and $\mathcal{%
D\subseteq M}_{m_{t}}\left( \mathbb{C}\right) .$ From the equation $\tau _{%
\mathcal{A}}\left( z\right) =\tau _{\mathcal{B}}\left( z\right) $ for all $%
z\in \mathcal{D},$ there is an integer $s$ such that $\mathcal{M}%
_{n_{t}}\left( \mathbb{C}\right) $ and $\mathcal{M}_{m_{t}}\left( \mathbb{C}%
\right) $ can be both unital embedded into $\mathcal{M}_{s}\left( \mathbb{C}%
\right) $ satisfying 
\begin{equation*}
\mathcal{M}_{s}\left( \mathbb{C}\right) \supseteq \mathcal{M}_{n_{t}}\left( 
\mathbb{C}\right) \supseteq \mathcal{D\subseteq M}_{m_{t}}\left( \mathbb{C}%
\right) \subseteq \mathcal{M}_{s}\left( \mathbb{C}\right) .
\end{equation*}%
Therefore, by Theorem \ref{70}, we have that 
\begin{equation*}
\left( \mathcal{M}_{n_{t}}\left( \mathbb{C}\right) ,E_{\mathcal{A}_{1}}|_{%
\mathcal{M}_{n_{t}}\left( \mathbb{C}\right) }\right) \ast _{\mathcal{D}%
}\left( \mathcal{M}_{m_{t}}\left( \mathbb{C}\right) ,E_{\mathcal{A}_{2}}|_{%
\mathcal{M}_{m_{t}}\left( \mathbb{C}\right) }\right) 
\end{equation*}%
is an MF algebra. 

Suppose $\left\{ x_{1},x_{2},\cdots ,z_{1},\cdots z_{h}\right\} $ is a
family of self-adjoint generators of $\mathcal{A}_{1}\ $with $\left\Vert
x_{i}\right\Vert =\left\Vert z_{j}\right\Vert =1$, $\left\{
y_{1},y_{2},\cdots z_{1},\cdots ,z_{h}\right\} $ is a family of self-adjoint
generators of $\mathcal{A}_{2}$ with $\left\Vert y_{i}\right\Vert
=\left\Vert z_{j}\right\Vert =1$ for all $i\in \mathbb{N}$ and $j\in \left\{
1,\cdots ,h\right\} $ where $\left\{ z_{1},\cdots ,z_{h}\right\} $ is a
family of self-adjoint generators of $\mathcal{D}.$ For any $\varepsilon >0$
and any $r\in \mathbb{N}$ with a family of noncommutative polynomials $%
P_{1},\cdots ,P_{r}$ in 
\begin{equation*}
\mathbb{C}\left\langle X_{1},\cdots X_{r},Y_{1},\cdots Y_{r},Z_{1},\cdots
,Z_{h}\right\rangle ,
\end{equation*}%
there is a integer $t\geq t_{0}$ with 
\begin{equation*}
T_{1},\cdots ,T_{r}\in \mathcal{M}_{m_{t}}^{s.a.}\left( \mathbb{C}\right)
\subseteq \left( \mathcal{M}_{n_{t}}\left( \mathbb{C}\right) ,E_{\mathcal{A}%
_{1}}|_{\mathcal{M}_{n_{t}}\left( \mathbb{C}\right) }\right) \ast _{\mathcal{%
D}}\left( \mathcal{M}_{m_{t}}\left( \mathbb{C}\right) ,E_{\mathcal{A}_{2}}|_{%
\mathcal{M}_{m_{t}}\left( \mathbb{C}\right) }\right) ,
\end{equation*}%
\begin{equation*}
S_{1},\cdots ,S_{r}\in \mathcal{M}_{n_{t}}^{s.a.}\left( \mathbb{C}\right)
\subseteq \left( \mathcal{M}_{n_{t}}\left( \mathbb{C}\right) ,E_{\mathcal{A}%
_{1}}|_{\mathcal{M}_{n_{t}}\left( \mathbb{C}\right) }\right) \ast _{\mathcal{%
D}}\left( \mathcal{M}_{m_{t}}\left( \mathbb{C}\right) ,E_{\mathcal{A}_{2}}|_{%
\mathcal{M}_{m_{t}}\left( \mathbb{C}\right) }\right) 
\end{equation*}%
and 
\begin{equation*}
K_{1},\cdots ,K_{h}\in \mathcal{D\subseteq }\left( \mathcal{M}_{n_{t}}\left( 
\mathbb{C}\right) ,E_{\mathcal{A}_{1}}|_{\mathcal{M}_{n_{t}}\left( \mathbb{C}%
\right) }\right) \ast _{\mathcal{D}}\left( \mathcal{M}_{m_{t}}\left( \mathbb{%
C}\right) ,E_{\mathcal{A}_{2}}|_{\mathcal{M}_{m_{t}}\left( \mathbb{C}\right)
}\right) 
\end{equation*}%
such that,$\ $for $1\leq j\leq r,$%
\begin{eqnarray*}
&&\left\vert \left\Vert P_{j}\left( x_{1},\cdots ,x_{r},y_{1},\cdots
,y_{r},z_{1},\cdots ,z_{h}\right) \right\Vert -\left\Vert P_{j}\left(
T_{1},\cdots ,T_{r},S_{1},\cdots ,S_{r},K_{1},\cdots ,K_{h}\right)
\right\Vert \right\vert  \\
&\leq &\left\Vert P_{j}\left( x_{1},\cdots ,x_{r},y_{1},\cdots
,y_{r},z_{1},\cdots ,z_{h}\right) -P_{j}\left( T_{1},\cdots
,T_{r},S_{1},\cdots ,S_{r},K_{1},\cdots ,K_{h}\right) \right\Vert  \\
&<&\varepsilon /2
\end{eqnarray*}%
where 
\begin{equation*}
P_{j}\left( T_{1},\cdots ,T_{r},S_{1},\cdots ,S_{r},K_{1},\cdots
,K_{h}\right) \in \left( \mathcal{M}_{n_{t}}\left( \mathbb{C}\right) ,E_{%
\mathcal{A}_{1}}|_{\mathcal{M}_{n_{t}}\left( \mathbb{C}\right) }\right) \ast
_{\mathcal{D}}\left( \mathcal{M}_{m_{t}}\left( \mathbb{C}\right) ,E_{%
\mathcal{A}_{2}}|_{\mathcal{M}_{m_{t}}\left( \mathbb{C}\right) }\right) .
\end{equation*}%
Since $\left( \mathcal{M}_{n_{t}}\left( \mathbb{C}\right) ,E_{\mathcal{A}%
_{1}}|_{\mathcal{M}_{n_{t}}\left( \mathbb{C}\right) }\right) \ast _{\mathcal{%
D}}\left( \mathcal{M}_{m_{t}}\left( \mathbb{C}\right) ,E_{\mathcal{A}_{2}}|_{%
\mathcal{M}_{m_{t}}\left( \mathbb{C}\right) }\right) $ is MF, there is a
positive integer $k$ and self-adjoint matrices $A_{1},\cdots
,A_{r},B_{1},\cdots ,B_{r},C_{1},\cdots ,C_{t}$ $\in \mathcal{M}%
_{k}^{s.a.}\left( \mathbb{C}\right) $ such that, for all $1\leq j\leq r$, 
\begin{equation*}
\left\vert \left\Vert P_{j}\left( A_{1},\cdots ,A_{r},B_{1},\cdots
,B_{r},C_{1},\cdots ,C_{h}\right) \right\Vert -\left\Vert P_{j}\left(
T_{1},\cdots ,T_{r},S_{1},\cdots ,S_{r},K_{1},\cdots ,K_{h}\right)
\right\Vert \right\vert \leq \varepsilon /2.
\end{equation*}%
This implies that 
\begin{equation*}
\left\vert \left\Vert P_{j}\left( x_{1},\cdots ,x_{r},y_{1},\cdots
,y_{r},z_{1},\cdots ,z_{h}\right) \right\Vert -\left\Vert P_{j}\left(
A_{1},\cdots ,A_{r},B_{1},\cdots ,B_{r},C_{1},\cdots ,C_{h}\right)
\right\Vert \right\vert <\varepsilon .
\end{equation*}%
Then, by Lemma 2.4.1 in \cite{DS2}, $\left( \mathcal{A}_{1},E_{\mathcal{A}%
_{1}}\right) \ast _{\mathcal{D}}\left( \mathcal{A}_{2},E_{\mathcal{A}%
_{2}}\right) $ is MF.
\end{proof}

We end this subsection with a result on the reduced C*-algebra of an
amalgamated free product group of two finite groups. First, we need the
following lemma. 

\begin{lemma}
\label{3.1}Suppose that $A\supseteq H\subseteq B$ are fnite groups
satisfying $\left[ A:H\right] \geq 3$ and $\left[ B:H\right] \geq 2.$ Let $%
A\ast _{H}B$ be the amalgamated free product of groups $A$ and $B$ over $H$.
Then $A\ast _{H}B$ is not amenable.
\end{lemma}

\begin{proof}
By Corollary 4.10 of \cite{BDJ}, we know that%
\begin{equation*}
\delta _{0}\left( \mathbb{C}\left( A\ast _{H}B\right) \right) =\left(
1-\left\vert A\right\vert ^{-1}\right) +\left( 1-\left\vert B\right\vert
^{-1}\right) -\left( 1-\left\vert H\right\vert ^{-1}\right) =1+\left(
\left\vert H\right\vert ^{-1}-\left\vert A\right\vert ^{-1}-\left\vert
B\right\vert ^{-1}\right) 
\end{equation*}%
where $\delta _{0}$ is the Voiculescu's (microstate) modifed free entropy
dimension and $\left\vert A\right\vert ,$ $\left\vert B\right\vert $ and $%
\left\vert H\right\vert $ are the cardinalities of groups $A,$ $B$ and $H$
respectively. Note that%
\begin{equation*}
\left\vert H\right\vert ^{-1}-\left\vert A\right\vert ^{-1}-\left\vert
B\right\vert ^{-1}=\left\vert H\right\vert ^{-1}\left( 1-\frac{1}{\left[ A:H%
\right] }-\frac{1}{\left[ B:H\right] }\right) \geq \left\vert H\right\vert
^{-1}\left( 1-\frac{1}{3}-\frac{1}{2}\right) >0.
\end{equation*}%
Hence $\delta _{0}\left( \mathbb{C}\left[ A\ast _{H}B\right] \right) >1.$ By
Corollary 6.8 in \cite{J}, we know that $A\ast _{H}B$ is not an amenable
group. 
\end{proof}

Applying preceding Theorem and Lemma, we obtain the following corollary.

\begin{corollary}
Suppose that $A\supseteq H\subseteq B$ are fnite groups and $A\ast _{H}B$ is
the amalgamated free product of groups $A$ and $B$ over $H$. Then $%
C_{r}^{\ast }\left( A\ast _{H}B\right) $ is an MF algebra. If futher $\left[
A:H\right] \geq 3$ and $\left[ B:H\right] \geq 2,$ then $Ext\left(
C_{r}^{\ast }\left( A\ast _{H}B\right) \right) ,$ the
Billmore-Douglus-Fillmore extension semigroup, is not a group.
\end{corollary}

\begin{proof}
Note that 
\begin{equation*}
C_{r}^{\ast }\left( A\ast _{H}B\right) \cong \left( C_{r}^{\ast }\left(
A\right) ,E_{A}\right) \ast _{C_{r}^{\ast }\left( H\right) }\left(
C_{r}^{\ast }\left( B\right) ,E_{B}\right) 
\end{equation*}%
where $E_{A}$ and $E_{B}$ are canonical trace preserving conditional
expectations from $C_{r}^{\ast }\left( A\right) $ onto $C_{r}^{\ast }\left(
H\right) ,$ and from $C_{r}^{\ast }\left( B\right) $ onto $C_{r}^{\ast
}\left( H\right) $ respectively (see \cite{BO}). Now it follows from Theorem %
\ref{70.1} that $C_{r}^{\ast }\left( A\ast _{H}B\right) $ is an MF algebra.
Moreover, if $\left[ A:H\right] \geq 3$ and $\left[ B:H\right] \geq 2,$ by
Lemma \ref{3.1}, we get that $A\ast _{H}B$ is not amenable, whence $%
C_{r}^{\ast }\left( A\ast _{H}B\right) $ is not a quasidiagonal C* -algebra.
Recall a well-known fact that if a unital C -algebra $\mathcal{E}$ is an MF
algebra but not a quasidiagonal C* -algebra, then Ext($\mathcal{E}$) is not
a group. Therefore, Ext($C_{r}^{\ast }\left( A\ast _{H}B\right) $), the
Billmore-Douglus-Fillmore extension semigroup, is not a group. 

{\small \addcontentsline{toc}{chapter}{\bf
BIBLIOGRAPHY}}
\end{proof}

Junhao Shen

Department of Mathematics and Statistics

Univerisity of New Hampshire

Durham, NH, 03824

Email: jog2@cisunix.unh.edu\bigskip 

Qihui Li

Department of Mathematics and Statistics

Univerisity of New Hampshire

Durham, NH, 03824

Email: qme2@cisunix.unh.edu\bigskip 

and

Department of Mathematics

East China University of Science and Technology

Shanghai, China Meilong Road 130, 200237

\bigskip 

\end{document}